\documentclass[12pt]{amsart}

\usepackage{amssymb}
\usepackage{graphicx}
\usepackage{ifthen}
\usepackage{pict2e}
\usepackage{xargs}
\usepackage{xspace}

\newcommand{\definedterm}[1]{\emph{#1}}

\newcommand{\absolutevalue}[1]{|#1|}

\newcommandx{\average}[2][2=]{\overline{#1}^{#2}}
\newcommand{\Bairespace}[1][]{
  \ifthenelse{\equal{#1}{}}{\functions{\N}{\N}}{\functions{#1}{\N}}
}
\newcommand{\Bairetree}[1][]{
  \ifthenelse{\equal{#1}{}}{\functions{<\N}{\N}}{\functions{#1}{\N}}
}

\newcommand{\calA}{\mathcal{A}}

\newcommand{\calN}{\mathcal{N}}

\newcommand{\calS}{\mathcal{S}}
\newcommand{\calU}{\mathcal{U}}

\newcommand{\Cantorspace}[1][]{
  \ifthenelse{\equal{#1}{}}{\functions{\N}{2}}{\functions{#1}{2}}
}
\newcommand{\Cantortree}[1][]{
  \ifthenelse{\equal{#1}{}}{\functions{<\N}{2}}{\functions{#1}{2}}
}
\newcommandx{\cardinality}[3][1 =, 2 =]{
  \ifthenelse{\equal{#1}{}}{|#3|}{|#3|_{#1}^{#2}}
}

\newcommand{\closedopeninterval}[2]{[#1, #2)}

\newcommand{\completion}[1]{\overline{#1}}
\newcommand{\composition}{\circ}

\newcommandx{\concatenation}[2][1 = undefined, 2 = undefined]{
  \ifthenelse{\equal{#1}{undefined}}{{}\smallfrown}{
    \ifthenelse{\equal{#2}{undefined}}{\bigoplus #1}{\bigoplus_{#1} #2}
  }
}
\newcommandx{\constant}[2][2 =]{\ifthenelse{\equal{#2}{}}{c_{#1}}{c_{#1, #2}}}
\newcommandx{\constantfunction}[3][2 =, 3 =]{
  \ifthenelse{\equal{#2}{}}{c \from #1 \to \image{c}{#1}}{c_{#3} \from #1 \to #2}
}

\newcommandx{\convolution}[2][1 = undefined, 2 = undefined]{
  \ifthenelse{\equal{#1}{undefined}}{\mathrel{*}}{
    \ifthenelse{\equal{#2}{undefined}}{\bigotimes #1}{\bigotimes_{#1} #2}
  }
}
\newcommandx{\Deltaclass}[2][1=,2=]{
  \ifthenelse{\equal{#2}{}}{\mathbf{\Delta}_{#1}}{\mathbf{\Delta}^{#1}_{#2}}
}

\newcommandx{\disjointunion}[2][1 =, 2 =]{
  \ifthenelse{\equal{#1}{}}{\sqcup}{
    \ifthenelse{\equal{#2}{}}{\bigsqcup #1}{{\bigsqcup_{#1} #2}}
  }
}
\newcommand{\displayedabsolutevalue}[1]{\left| #1 \right|}

\newcommand{\equivalenceclass}[2]{[#1]_{#2}}

\newcommand{\extendedextensions}[1]{\extensions{#1}^*}

\newcommand{\extension}[1]{\overline{#1}}
\newcommand{\extensions}[2][]{
  \ifthenelse{\equal{#1}{}}{\calN_{#2}}{\extendedextensions{#2} \intersection #1}
}

\newcommand{\from}{\colon}

\newcommandx{\functions}[3][3 =]{
  \ifthenelse{\equal{#3}{}}{#2^{#1}}{#2^{#1}_{#3}}
}

\newcommand{\goesto}{\rightarrow}
\newcommand{\graph}[1]{\mathrm{graph}(#1)}

\newcommand{\GzeroN}[1][]{
  \ifthenelse{\equal{#1}{}}{\mathbb{G}_0^\N}{\mathbb{G}_{0, #1}^{\N}}
}

\newcommand{\identity}[1]{1_{#1}}
\newcommand{\image}[2]{#1(#2)}

\newcommandx{\identityfunction}[2][2 =]{
  \ifthenelse{\equal{#2}{}}{\mathrm{id} \from #1 \to #1}{\mathrm{id} \from #1 \to #2}
}
\newcommand{\infimum}[2][]{
  \ifthenelse{\equal{#1}{}}{\inf #1}{\inf_{#1}{#2}}
}
\newcommand{\injections}[1]{\functions{\N}{(\N)}}
\newcommandx{\intersection}[2][1 =, 2 =]{
  \ifthenelse{\equal{#1}{}}{\cap}{
    \ifthenelse{\equal{#2}{}}{\bigcap #1}{{\bigcap_{#1} #2}}
  }
}
\newcommand{\inverse}[1]{#1^{-1}}

\newcommand{\lacunaritygraph}[2]{G^{#2}_{#1}}
\newcommand{\length}[1]{|#1|}
\newcommandx{\limit}[2][1 =, 2 =]{
  \ifthenelse{\equal{#1}{}}{\lim}{
    \ifthenelse{\equal{#2}{}}{\lim #1}{{\lim_{#1} #2}}
  }
}

\newcommand{\measure}[2]{\mu_{#1}^{#2}}

\newcommand{\N}{\mathbb{N}}

\newcommand{\openinterval}[2]{(#1, #2)}

\newcommand{\orbitequivalencerelation}[2]{E_{#1}^{#2}}
\newcommand{\pair}[2]{(#1, #2)}
\newcommandx{\Piclass}[2][1=,2=]{
  \ifthenelse{\equal{#2}{}}{\mathbf{\Pi}_{#1}}{\mathbf{\Pi}^{#1}_{#2}}
}

\newcommand{\preimage}[2]{#1^{-1}(#2)}

\newcommandx{\product}[2][1 =, 2 =]{
  \ifthenelse{\equal{#1}{}}{\times}{
    \ifthenelse{\equal{#2}{}}{\prod #1}{{\prod_{#1} #2}}
  }
}
\newcommandx{\projection}[1]{\mathrm{proj}_{#1}}
\newcommand{\pushforward}[2]{#1_* #2}

\renewcommandx{\restriction}[3][3 = undefined]{
  \ifthenelse{\equal{#3}{undefined}}{#1 \upharpoonright #2}{#1
    \upharpoonright #2 \to #3}
}
\newcommand{\saturation}[2]{[#1]_{#2}}
\newcommandx{\sequence}[2][2 = undefined]{
  \ifthenelse{\equal{#2}{undefined}}{(#1)}{
    (#1)_{#2}
  }
}
\newcommandx{\set}[2][2 = undefined]{
  \ifthenelse{\equal{#2}{undefined}}{\{ #1 \}}{
    \{ #1 \suchthat #2 \}
  }
}
\newcommand{\setcomplement}[1]{\twiddle #1}
\newcommandx{\sets}[3][3 =]{
  \ifthenelse{\equal{#3}{}}{[#2]^{#1}}{[#2]^{#1}_{#3}}
}

\newcommandx{\Sigmaclass}[2][1=,2=]{
  \ifthenelse{\equal{#2}{}}{\mathbf{\Sigma}_{#1}}{\mathbf{\Sigma}^{#1}_{#2}}
}

\newcommand{\suchthat}{\mid}

\newcommand{\twiddle}{\raisebox{1pt}{\scalebox{.75}{$\mathord{\sim}$}}}
\newcommandx{\union}[2][1 =, 2 =]{
  \ifthenelse{\equal{#1}{}}{\cup}{
    \ifthenelse{\equal{#2}{}}{\bigcup #1}{{\bigcup_{#1} #2}}
  }
}
\newcommand{\verticalsection}[2]{#1_{#2}}

\newcommand{\Z}{\mathbb{Z}}

\newcommand{\Becker}{Beck\-er\xspace}
\newcommand{\Borel}{Bor\-el\xspace}

\newcommand{\Caratheodory}{Car\-ath\-\'{e}o\-dory\xspace}
\newcommand{\Dougherty}{Dough\-er\-ty\xspace}
\newcommand{\Effros}{Eff\-ros\xspace}
\newcommand{\Feldman}{Feld\-man\xspace}
\newcommand{\Glimm}{Glimm\xspace}

\newcommand{\Hopf}{Hopf\xspace}

\newcommand{\Jackson}{Jack\-son\xspace}

\newcommand{\Kechris}{Kech\-ris\xspace}
\newcommand{\Koenig}{K\"{o}nig\xspace}

\newcommand{\Lusin}{Lu\-sin\xspace}

\newcommand{\Moore}{Moore\xspace}
\newcommand{\Murray}{Mur\-ray\xspace}
\newcommand{\Nadkarni}{Nad\-kar\-ni\xspace}
\newcommand{\Nikodym}{Nik\-o\-dym\xspace}
\newcommand{\Novikov}{No\-vik\-ov\xspace}
\newcommand{\Polish}{Po\-lish\xspace}
\newcommand{\Radon}{Ra\-don\xspace}

\newcommand{\vonNeumann}{von Neu\-mann\xspace}

\newenvironment{lemmaproof}{
  
  \begin{proof}
}{\end{proof}}

\newenvironment{propositionproof}{
  
  \begin{proof}
}{\end{proof}}

\newenvironment{sublemmaproof}{
  
  \begin{proof}
}{\end{proof}}

\newenvironment{theoremproof}{
  
  \begin{proof}
}{\end{proof}}

\newtheorem{lemma}{Lemma}[section]
\newtheorem{proposition}[lemma]{Proposition}
\newtheorem{sublemma}[lemma]{Sublemma}
\newtheorem{theorem}[lemma]{Theorem}

\newtheorem{introtheorem}{Theorem}

\theoremstyle{definition}
\newtheorem{remark}[lemma]{Remark}

\newenvironment{acknowledgements}{
  \textbf{Acknowledgements.}
}

\begin{document}


\begin{abstract}
  We show that a natural generalization of compressibility 
  is the sole obstruction to the existence of a cocycle-invariant
  \Borel probability measure.
\end{abstract}

\author[B.D. Miller]{Benjamin D. Miller}

\address{
  Benjamin D. Miller \\
  Kurt G\"{o}del Research Center for Mathematical Logic \\
  Universit\"{a}t Wien \\
  W\"{a}hringer Stra{\ss}e 25 \\
  1090 Wien \\
  Austria
 }

\email{benjamin.miller@univie.ac.at}

\urladdr{
  http://www.logic.univie.ac.at/benjamin.miller
}

\thanks{The author was supported in part by FWF Grants
  P28153 and P29999.}
  


\title[Cocycle-invariant probability measures]{On the existence of
  cocycle-invariant Borel probability measures}

\maketitle

\section*{Introduction}

Suppose that $X$ is a standard \Borel space and $T \from X \to X$
is a \Borel automorphism of $X$. A Borel measure $\mu$ on $X$ is
\definedterm{$T$-invariant} if $\mu(\image{T}{B}) = \mu(B)$ for all
\Borel sets $B \subseteq X$. The characterization of the class of
\Borel automorphisms of standard \Borel spaces admitting an invariant
\Borel probability measure is a fundamental problem going back to
\Hopf (see \cite{Hopf}). 

A \definedterm{compression} of an equivalence relation $E$ on $X$
is an injection $\phi \from X \to X$ sending each $E$-class into a
proper subset of itself. Building on work of \Murray-\vonNeumann
(see \cite{MurrayVonNeumann}), \Nadkarni has shown that the
existence of a \Borel compression of the orbit equivalence
relation $\orbitequivalencerelation{T}{X}$ induced by $T$ is the sole
obstruction to the existence of a $T$-invariant \Borel probability
measure (see \cite{Nadkarni}).

Suppose that $E$ is a \Borel equivalence relation on $X$ that is
\definedterm{countable}, in the sense that all of its equivalence
classes are countable. A \Borel measure $\mu$ on $X$ is \definedterm
{$E$-invariant} if it is $T$-invariant for all \Borel automorphisms
$T \from X \to X$ whose graphs are contained in $E$. It is easy to
see that a \Borel measure is $T$-invariant if and only if it is
$\orbitequivalencerelation{T}{X}$-invariant. \Becker-\Kechris have
pointed out that \Nadkarni's argument yields the more general fact
that the existence of a \Borel compression of $E$ is the sole
obstruction to the existence of an $E$-invariant \Borel probability
measure (see \cite[Theorem 4.3.1]{BeckerKechris}).

An equivalence relation is \definedterm{aperiodic} if all of its classes
are infinite. A set $Y \subseteq X$ is \definedterm{$E$-complete} if
it intersects every $E$-class in at least one point, and a set $Y
\subseteq X$ is a \definedterm{partial transversal} of $E$ if it
intersects every $E$-class in at most one point. A \definedterm
{transversal} of $E$ is an $E$-complete partial transversal of $E$.
The \Lusin-\Novikov uniformization theorem (see, for example,
\cite[Theorem 18.10]{Kechris}) ensures that there is a \Borel
transversal of $E$ if and only if $X$ is the union of countably-many
\Borel partial transversals of $E$. We say that $E$ is \definedterm
{smooth} if it satisfies these equivalent conditions.
\Dougherty-\Jackson-\Kechris have pointed out that the
existence of a \Borel compression of $E$ is equivalent to the
existence of an aperiodic smooth \Borel subequivalence relation of
$E$ (see \cite[Proposition 2.5]{DoughertyJacksonKechris}), thereby
obtaining another characterization of the class of countable \Borel
equivalence relations on standard \Borel spaces admitting an
invariant \Borel probability measure.

A substantially weaker notion than $E$-invariance is that of
\definedterm{$E$-quasi-invariance}, where one asks that
$\mu(\image{T}{B}) = 0 \iff \mu(B) = 0$ for all \Borel sets $B
\subseteq X$ and \Borel automorphisms $T \from X \to X$ whose
graphs are contained in $E$. Given a group $\Gamma$, we say that
a function $\rho \from E \to \Gamma$ is a \definedterm{cocycle} if
$\rho(x, z) = \rho(x, y) \rho(y, z)$ whenever $x \mathrel{E} y \mathrel
{E} z$. Given a \Borel cocycle $\rho \from E \to \openinterval{0}
{\infty}$, we say that a \Borel measure $\mu$ on $X$ is \definedterm
{$\rho$-invariant} if $\mu(\image{T}{B}) = \int_B \rho(T(x), x)
\ d\mu(x)$ for all \Borel sets $B \subseteq X$ and \Borel
automorphisms $T \from X \to X$ whose graphs are contained in
$E$. Clearly $E$-invariance is equivalent to invariance with respect
to the constant cocycle, whereas the \Radon-\Nikodym Theorem
(see, for example, \cite[{\S}17.A]{Kechris}) and the \Feldman-\Moore
observation that countable \Borel equivalence relations on standard
\Borel spaces are orbit equivalence relations induced by \Borel
actions of countable groups (see \cite[Theorem 1]
{FeldmanMoore}) ensure that $E$-quasi-invariance is equivalent to
invariance with respect to some \Borel cocycle $\rho \from E \to
\openinterval{0}{\infty}$ (see, for example, \cite[\S8]{KechrisMiller}).
A characterization of the class of \Borel cocycles $\rho \from E \to
\openinterval{0}{\infty}$ admitting an invariant \Borel probability
measure was provided in \cite{Miller:Probability}. Here we
investigate more natural generalizations of the characterizations
mentioned above.

In \S\ref{smooth}, we introduce the direct generalizations of
\definedterm{aperiodicity} and \definedterm{compressibility} to
cocycles that come from viewing $\rho$ as endowing each
$E$-class with a notion of relative size. We also introduce the
generalization of \definedterm{smoothness} to cocycles that comes
from the \Glimm-\Effros dichotomy. We note that, unfortunately, even
when $E$ is smooth, there are \Borel cocycles on $E$ admitting
neither a compression nor an invariant \Borel probability measure. In
order to bypass this obstacle, we introduce the \definedterm{quotient}
of $\rho$ by a finite subequivalence relation of $E$. Generalizing the
observation of \Dougherty-\Jackson-\Kechris, we show that the
existence of an injective \Borel compression of the quotient of $\rho$
by a finite \Borel subequivalence relation of $E$ is equivalent to the
existence of a \Borel subequivalence relation of $E$ on which $\rho$
is aperiodic and smooth. We also note that, at least when $\rho$ is
smooth, the existence of an injective \Borel compression of the
quotient of $\rho$ by a finite \Borel subequivalence relation of $E$ is
the sole obstacle to the existence of a $\rho$-invariant \Borel
probability measure.

In \S\ref{coboundaries}, we introduce \definedterm{\Borel
coboundaries}, a natural class of particularly simple \Borel cocycles
containing the constant cocycles. We note that, unfortunately, there
are \Borel coboundaries admitting neither an injective \Borel
compression of the quotient by a finite \Borel subequivalence relation
of $E$ nor an invariant \Borel probability measure. In order to bypass
this new obstacle, we then drop the assumption of injectivity, and
combine the \Becker-\Kechris generalization of \Nadkarni's theorem,
the \Dougherty-\Jackson-\Kechris characterization of the existence
of \Borel compressions, and an approximation lemma to generalize
\Nadkarni's theorem to \Borel coboundaries.

\begin{introtheorem} \label{introduction:coboundaries}
  Suppose that $X$ is a standard \Borel space, $E$ is a countable
  \Borel equivalence relation on $X$, and $\rho \from E \to
  \openinterval{0}{\infty}$ is a \Borel coboundary. Then exactly one
  of the following holds:
  \begin{enumerate}
    \item There is a finite-to-one \Borel compression of the quotient of
      $\rho$ by a finite \Borel subequivalence relation of $E$.
    \item There is a $\rho$-invariant \Borel probability measure.
  \end{enumerate}
\end{introtheorem}

In \S\ref{cocycles}, we no longer restrict our attention to \Borel
coboundaries. Unfortunately, the direct generalization of Theorem
\ref{introduction:coboundaries} to \Borel cocycles remains open. In
order to bypass this final obstacle, we consider the weakening of the
notion of a compression of the quotient of $\rho$ by a finite
subequivalence relation $F$ of $E$ obtained by only taking the
quotient in the range, which we refer to as a \definedterm
{compression} of $\rho$ over $F$. By augmenting the main argument
of \cite{Miller:Probability} with an additional approximation lemma, we
generalize \Nadkarni's theorem to \Borel cocycles.

\begin{introtheorem}
  Suppose that $X$ is a standard \Borel space, $E$ is a countable
  \Borel equivalence relation on $X$, and $\rho \from E \to
  \openinterval{0}{\infty}$ is a \Borel cocycle. Then exactly one
  of the following holds:
  \begin{enumerate}
    \item There is a finite-to-one \Borel compression of $\rho$ over a
      finite \Borel subequivalence relation of $E$.
    \item There is a $\rho$-invariant \Borel probability measure.
  \end{enumerate}
\end{introtheorem}

\section{Smooth cocycles} \label{smooth}

One can think of a cocycle $\rho \from E \to \openinterval{0}{\infty}$
as assigning a notion of relative size to each $E$-class $C$, with the
\definedterm{$\rho$-size} of a point $y \in C$ relative to a point $z \in
C$ being $\rho(y, z)$. More generally, the \definedterm{$\rho$-size}
of a set $Y \subseteq C$ relative to $z$ is given by $\cardinality[z]
[\rho]{Y} = \sum_{y \in Y} \rho(y, z)$. We say that $Y$ is \definedterm
{$\rho$-infinite} if this quantity is infinite. As the definition of cocycle
ensures that $\cardinality[z'][\rho]{Y} = \cardinality[z][\rho]{Y}
\thinspace \rho(z, z')$ for all $z' \in C$, it follows that the notion of
being $\rho$-infinite does not depend on the choice of $z \in C$. It
also follows that the \definedterm{$\rho$-size} of $Y$ relative to a
non-empty set $Z \subseteq C$, given by $\cardinality[Z][\rho]{Y} =
\cardinality[z][\rho]{Y} / \cardinality[z][\rho]{Z}$, does not depend on
the choice of $z \in C$.

We say that a cocycle $\rho \from E \to \openinterval{0}{\infty}$ is
\definedterm{aperiodic} if every $E$-class is $\rho$-infinite. Note
that the aperiodicity of $\rho$ trivially yields that of $E$. Conversely,
when $\rho$ is bounded, the aperiodicity of $E$ yields that of $\rho$.

We say that a function $\phi \from X \to X$ is a \definedterm
{compression} of $\rho$ if the graph of $\phi$ is contained in $E$,
$\cardinality[x][\rho]{\preimage{\phi}{x}} \le 1$ for all $x \in X$, and
the set $\set{x \in X}[{\cardinality[x][\rho]{\preimage{\phi}{x}} < 1}]$
is $E$-complete. Note that, when $\rho$ is the constant cocycle,
a function $\phi \from X \to X$ is a compression of $E$ if and only
if it is a compression of $\rho$.

\begin{proposition} \label{smooth:counterexample}
  Suppose that $X$ is a standard \Borel space and $E$ is an
  aperiodic smooth countable \Borel equivalence relation on $X$.
  Then there is an aperiodic \Borel cocycle $\rho \from E \to
  \openinterval{0}{\infty}$ that does not admit a compression.
\end{proposition}

\begin{propositionproof}
  Fix a strictly decreasing sequence $\sequence{r_n}[n \in \N]$ of
  positive real numbers for which $\sum_{n \in \N} r_n = \infty$. As
  $E$ is both aperiodic and smooth, the \Lusin-\Novikov
  uniformization theorem yields a partition $\sequence{B_n}[n \in \N]$
  of $X$ into \Borel transversals of $E$. For each $x \in X$, let $n(x)$
  denote the unique natural number for which $x \in B_{n(x)}$, and
  define $\rho \from E \to \openinterval{0}{\infty}$ by setting $\rho(x,
  y) = r_{n(x)} / r_{n(y)}$ whenever $x \mathrel{E} y$. 
  
  The fact that $\sum_{n \in \N} r_n = \infty$ ensures that $\rho$ is
  aperiodic. To see that there is no compression of $\rho$, note that
  if $\phi \from X \to X$ is a function such that the graph of $\phi$ is
  contained in $E$ and $\cardinality[x][\rho]{\preimage{\phi}{x}} \le 1$
  for all $x \in X$, then a straightforward induction on $n(x)$, using
  the fact that $\sequence{r_n}[n \in \N]$ is strictly decreasing, shows
  that $\phi(x) = x$ for all $x \in X$.
\end{propositionproof}

A \definedterm{digraph} on $X$ is an irreflexive set $G \subseteq X
\times X$. Given such a digraph, we say that a set $Y \subseteq X$
is \definedterm{$G$-independent} if $G \intersection (Y \times Y) =
\emptyset$. A \definedterm{$Y$-coloring} of $G$ is a function $c
\from X \to Y$ with the property that $\preimage{c}{y}$ is
$G$-independent for all $y \in Y$.

The \definedterm{vertical sections} of a set $R \subseteq X \times
Y$ are the sets of the form $\verticalsection{R}{x} = \set{y \in Y}
[\pair{x}{y} \in R]$, where $x \in X$. When $G$ is \Borel, it follows
from \cite[Proposition 4.5]{KechrisSoleckiTodorcevic} that there is
a \Borel $\N$-coloring of $G$ if and only if $X$ is the union of
countably-many \Borel sets $B \subseteq X$ for which the vertical
sections of $G \intersection (B \times B)$ are finite.

We say that a \Borel measure $\mu$ on $X$ is \definedterm
{$E$-ergodic} if every $E$-invariant \Borel set is $\mu$-conull or
$\mu$-null. Given a \Borel cocycle $\rho \from E \to \Gamma$ and a
set $Z \subseteq \Gamma$, let $\lacunaritygraph{Z}{\rho}$ denote
the digraph on $X$ with respect to which distinct points $x$ and
$y$ are related if and only if they are $E$-equivalent and $\rho(x, y) \in
Z$. The \Glimm-\Effros dichotomy for countable \Borel equivalence
relations (see \cite{Weiss}) ensures that $E$ is smooth if and only if
there is no atomless $E$-ergodic $E$-invariant $\sigma$-finite
\Borel measure. In \cite{Miller:Sigma}, this was generalized to show
that if $\rho \from E \to \openinterval{0}{\infty}$ is a \Borel cocycle,
then there is an open neighborhood $U \subseteq \openinterval{0}
{\infty}$ of $1$ for which there is a \Borel $\N$-coloring of
$\lacunaritygraph{U}{\rho}$ if and only if there is no atomless
$E$-ergodic $\rho$-invariant $\sigma$-finite \Borel measure.
Consequently, we say that a \Borel cocycle $\rho \from E \to
\openinterval{0}{\infty}$ is \definedterm{smooth} if it satisfies these
equivalent conditions. Note that the smoothness of $E$ trivially
yields that of $\rho$. Conversely, when $\rho$ is bounded, the
smoothness of $\rho$ ensures that $X$ is the union of
countably-many \Borel sets whose intersection with each $E$-class
is finite, thus $E$ is smooth.

We say that a set $Y \subseteq X$ is \definedterm{$\rho$-lacunary}
if it is $\lacunaritygraph{U}{\rho}$-independent for some open
neighborhood $U \subseteq \openinterval{0}{\infty}$ of $1$.

\begin{proposition} \label{smooth:openyieldscompact}
  Suppose that $X$ is a standard \Borel space, $E$ is a countable
  \Borel equivalence relation on $X$, $\Gamma$ is a \Polish group,
  and $\rho \from E \to \Gamma$ is a \Borel cocycle. If there is an
  open neighborhood $U \subseteq \Gamma$ of $\identity
  {\Gamma}$ for which there is a \Borel $\N$-coloring of
  $\lacunaritygraph{U}{\rho}$, then there is a \Borel $\N$-coloring of
  $\lacunaritygraph{K}{\rho}$ for all compact sets $K \subseteq
  \Gamma$.
\end{proposition}

\begin{propositionproof}
  Given a digraph $G$ on $X$, we say that a set $Y \subseteq X$ is
  a \definedterm{$G$-clique} if all pairs of distinct points of $Y$ are
  $G$-related. It is sufficient to show that if a set $Y \subseteq X$ does
  not contain an infinite $\lacunaritygraph{U}{\rho}$-clique, then the
  vertical sections of $\lacunaritygraph{K}{\rho} \intersection (X \times
  Y)$ are finite. Towards this end, fix a non-empty open set $V \subseteq
  \Gamma$ with the property that $\inverse{V} V \subseteq U$, as
  well as a finite sequence $\sequence{\gamma_i}[i < n]$ of elements
  of $\Gamma$ for which $K \subseteq \union[i < n][\gamma_i V]$,
  and note that if $x \in X$, then $\verticalsection{(\lacunaritygraph
  {K}{\rho})}{x} \subseteq \union[i < n][\verticalsection{(\lacunaritygraph
  {\gamma_i V}{\rho})}{x}]$, so we need only show that each
  $\verticalsection{(\lacunaritygraph{\gamma_i V}{\rho})}{x}$ is a
  $\lacunaritygraph{U}{\rho}$-clique. But if $i < n$ and $y, z \in
  \verticalsection{(\lacunaritygraph{\gamma_i V}{\rho})}{x}$, then
  $\rho(y, z) = \rho(y, x) \rho(x, z) \in \inverse{(\gamma_i V)} \gamma_i
  V = \inverse{V} V \subseteq U$.
\end{propositionproof}

The following fact ensures that a \Borel cocycle $\rho \from E
\to \openinterval{0}{\infty}$ is smooth if and only if there is an
$E$-complete $\rho$-lacunary \Borel set.

\begin{proposition} \label{smooth:completeyieldscoloring}
  Suppose that $X$ is a standard \Borel space, $E$ is a countable
  \Borel equivalence relation on $X$, $\Gamma$ is a locally compact
  \Polish group, $\rho \from E \to \Gamma$ is a \Borel cocycle, and
  $U \subseteq \Gamma$ is a pre-compact open neighborhood of
  $\identity{\Gamma}$. Then there is a \Borel $\N$-coloring of
  $\lacunaritygraph{U}{\rho}$ if and only if there is an $E$-complete
  $\lacunaritygraph{U}{\rho}$-independent \Borel set. 
\end{proposition}

\begin{propositionproof}
  If $c \from X \to \N$ is a \Borel $\N$-coloring of $\lacunaritygraph
  {U}{\rho}$, then set $A_n = \preimage{c}{n}$ and $B_n = A_n
  \setminus \union[m < n][\saturation{A_m}{E}]$ for all $n \in \N$.
  As the \Lusin-\Novikov uniformization theorem ensures that the
  latter sets are \Borel, it follows that their union is an $E$-complete
  $\lacunaritygraph{U}{\rho}$-independent \Borel set. 

  Conversely, suppose that $B \subseteq X$ is an $E$-complete
  $\lacunaritygraph{U}{\rho}$-independent \Borel set. The
  \Lusin-\Novikov uniformization theorem then yields \Borel functions
  $\phi_n \from B \to X$ such that $E \intersection (B \times X) =
  \union[n \in \N][\graph{\phi_n}]$, from which it follows that there are
  such functions satisfying the additional constraint that the sets
  $K_n = \image{\rho}{\graph{\phi_n}}$ are pre-compact. As
  Proposition \ref{smooth:openyieldscompact} yields \Borel
  $\N$-colorings of $\lacunaritygraph{K_n U \inverse{K_n}}{\rho}
  \intersection (B \times B)$, and the \Lusin-\Novikov uniformization
  theorem ensures that $\phi_n$ sends $\lacunaritygraph{K_n U
  \inverse{K_n}}{\rho}$-independent \Borel sets to $\lacunaritygraph
  {U}{\rho}$-independent \Borel sets, there are \Borel
  $\N$-colorings of $\lacunaritygraph{U}{\rho} \intersection
  (\image{\phi_n}{B} \times \image{\phi_n}{B})$, and therefore of
  $\lacunaritygraph{U}{\rho}$.
\end{propositionproof}

\begin{remark}
  Propositions \ref{smooth:openyieldscompact} and \ref
  {smooth:completeyieldscoloring} easily imply that a \Borel cocycle
  $\rho \from E \to \openinterval{0}{\infty}$ is smooth if and only if $X$
  is the union of countably-many $\rho$-lacunary \Borel sets.
\end{remark}

We say that a function $\phi \from X \to X$ is \definedterm{strictly
$\rho$-increasing} if its graph is contained in $E$ and $\cardinality
[x][\rho]{\preimage{\phi}{x}} < 1$ for all $x \in X$.

\begin{proposition} \label{smooth:smoothandstrictlyincreasing}
  Suppose that $X$ is a standard \Borel space, $E$ is a countable
  \Borel equivalence relation on $X$, and $\rho \from E \to
  \openinterval {0}{\infty}$ is a smooth \Borel cocycle. Then there
  is an $E$-invariant \Borel set $B \subseteq X$ for which
  $\restriction{E}{\setcomplement{B}}$ is smooth and there is a
  strictly $(\restriction{\rho}{(\restriction{E}{B})})$-increasing \Borel
  automorphism.
\end{proposition}

\begin{propositionproof}
  Fix a partition $\sequence{B_n}[n \in \N]$ of $X$ into
  $\rho$-lacunary \Borel sets. For each $x \in X$, let $n(x)$ be
  the unique natural number for which $x \in B_{n(x)}$. Let $\preceq$ be
  the partial order on $X$ with respect to which $x \preceq y$ if and
  only if $x \mathrel{E} y$, $n(x) = n(y)$, and $\rho(x, y) \le 1$, and
  let $B$ be the set of $x \in X$ such that for all $n \in \N$, either
  $B_n \intersection \equivalenceclass{x}{E} = \emptyset$ or
  $\restriction{\mathord{\preceq}}{(B_n \intersection \equivalenceclass
  {x}{E})}$ is isomorphic to the usual ordering of $\Z$. Then
  $\restriction{E}{\setcomplement{B}}$ is smooth, and the
  $(\restriction{\mathord{\preceq}}{B})$-successor function is a strictly
  $(\restriction{\rho}{(\restriction{E} {B})})$-increasing \Borel
  automorphism.
\end{propositionproof}

Given a cocycle $\rho \from E \to \openinterval{0}{\infty}$ and a
finite subequivalence relation $F$ of $E$, define $\rho / F \from E / F
\to \openinterval{0}{\infty}$ by $(\rho / F)(\equivalenceclass{x}{F},
\equivalenceclass{y}{F}) = \cardinality[\equivalenceclass{y}{F}][\rho]
{\equivalenceclass{x}{F}}$. The \Lusin-\Novikov uniformization theorem
ensures that if $F$ is \Borel, then $X / F$ is standard \Borel, so that
$E / F$ is a countable \Borel equivalence relation on a standard
\Borel space. Moreover, if $\rho$ is \Borel, then $\rho / F$ is a \Borel
cocycle on $E / F$. The \Lusin-\Novikov uniformization theorem also
implies that, when $\rho$ is the constant cocycle, a \Borel
compression of $\rho / F$ gives rise to a \Borel compression of
$\rho$. In spite of Proposition \ref{smooth:counterexample}, such
quotients allow us to generalize the fact that aperiodic smooth
countable \Borel equivalence relations admit \Borel compressions.

\begin{proposition} \label{smooth:aperiodic}
  Suppose that $X$ is a standard \Borel space, $E$ is a countable
  \Borel equivalence relation on $X$, and $\rho \from E \to
  \openinterval{0}{\infty}$ is an aperiodic smooth \Borel cocycle.
  Then there is a finite \Borel subequivalence relation $F$ of $E$
  for which there is a strictly $(\rho / F)$-increasing \Borel injection.
\end{proposition}

\begin{propositionproof}
  By Proposition \ref{smooth:smoothandstrictlyincreasing}, we can
  assume that $E$ is smooth. As the aperiodicity of $\rho$ yields
  that of $E$, there is a partition $\sequence{B_n}[n \in \N]$ of $X$
  into \Borel transversals of $E$. For each $x \in X$, let $n(x)$
  be the unique natural number with $x \in B_{n(x)}$, set
  $n_i(x) = i$ for all $i < 2$, recursively define $n_{i+2}(x)$ to
  be the least natural number such that the $\rho$-size of the set
  $\set{y \in \equivalenceclass{x}{E}}[n_{i+1}(x) \le n(y) < n_{i+2}(x)]$
  relative to the set $\set{y \in \equivalenceclass{x}{E}}[n_i(x) \le n(y)
  < n_{i+1}(x)]$ is strictly greater than one for all $i \in \N$, and let
  $i(x)$ be the unique natural number with the property that $n_{i(x)}
  (x) \le n(x) < n_{i(x) + 1}(x)$. Let $F$ be the subequivalence
  relation of $E$ with respect to which two $E$-equivalent points are
  $F$-equivalent if and only if $i(x) = i(y)$. Then the function $\phi
  \from X / F \to X / F$, given by $\phi(\equivalenceclass{x}{F}) = \set
  {y \in \equivalenceclass{x}{E}}[i(y) = i(x) + 1]$, is a strictly $(\rho /
  F)$-increasing \Borel injection.
\end{propositionproof}

The following fact yields an equivalent form of $\rho$-invariance
that will prove useful when considering \Borel injections.

\begin{proposition} \label{smooth:injection}
  Suppose that $X$ is a standard \Borel space, $E$ is a countable
  \Borel equivalence relation on $X$, $\rho \from E \to \openinterval
  {0}{\infty}$ is a \Borel cocycle, and $\mu$ is a $\rho$-invariant
  \Borel measure. Then $\mu(\image{T}{B}) = \int_B \rho(T(x), x)
  \ d\mu(x)$ for all \Borel sets $B \subseteq X$ and \Borel injections
  $T \from B \to X$ whose graphs are contained in $E$.
\end{proposition}

\begin{propositionproof}
  Fix a countable group $\Gamma = \set{\gamma_n}[n \in \N]$ of
  \Borel automorphisms of $X$ whose induced orbit equivalence
  relation is $E$, recursively define $B_n = \set{x \in B \setminus
  \union[m < n][B_m]}[T(x) = \gamma_n \cdot x]$ for all $n \in \N$,
  and note that
  \begin{align*}
    \mu(\image{T}{B})
      & = \sum_{n \in \N} \mu(\image{\gamma_n}{B_n}) \\
      & = \sum_{n \in \N} \int_{B_n} \rho(\gamma_n \cdot x, x)
        \ d\mu(x) \\
      & = \int_B \rho(T(x), x) \ d\mu(x)
  \end{align*}
  by $\rho$-invariance.
\end{propositionproof}

The following fact yields an equivalent form of $\rho$-invariance
that will prove useful when considering \Borel functions.

\begin{proposition} \label{smooth:invariance}
  Suppose that $X$ is a standard \Borel space, $E$ is a countable
  \Borel equivalence relation on $X$, $\rho \from E \to \openinterval
  {0}{\infty}$ is a \Borel cocycle, and $\mu$ is a $\rho$-invariant
  \Borel measure. Then $\mu(\preimage{\phi}{B}) = \int_B
  \cardinality[x][\rho]{\preimage{\phi}{x}} \ d\mu(x)$ for all \Borel sets
  $B \subseteq X$ and \Borel functions $\phi \from X \to X$ whose
  graphs are contained in $E$.
\end{proposition}

\begin{propositionproof}
  By the \Lusin-\Novikov uniformization theorem, there are \Borel
  sets $B_n \subseteq B$ and \Borel injections $T_n \from B_n \to
  X$ with the property that $\sequence{\graph{T_n}}[n \in \N]$
  partitions $\graph{\inverse{\phi}} \intersection (B \times X)$. Then
  \begin{equation*}
    \int_B \cardinality[x][\rho]{\preimage{\phi}{x}} \ d\mu(x)
      = \sum_{n \in \N} \int_{B_n} \rho(T_n(x), x)
        \ d\mu(x)
      = \mu(\preimage{\phi}{B})
  \end{equation*}
  by Proposition \ref{smooth:injection}.
\end{propositionproof}

Much as before, we say that a function $\phi \from X \to X$ is a
\definedterm{compression} of $\rho$ over a finite subequivalence
relation $F$ of $E$ if the graph of $\phi$ is contained in $E$,
$\cardinality[\equivalenceclass{x}{F}][\rho]{\preimage{\phi}
{\equivalenceclass{x}{F}}} \le 1$ for all $x \in X$, and the set $\set{x
\in X}[{\cardinality[\equivalenceclass{x}{F}][\rho]{\preimage{\phi}
{\equivalenceclass{x}{F}}} < 1}]$ is $E$-complete. The
\Lusin-\Novikov uniformization theorem ensures that every \Borel
compression of the quotient of $\rho$ by a finite \Borel
subequivalence relation $F$ of $E$ gives rise to a \Borel
compression of $\rho$ over $F$. It also implies that, when $\rho$ is
the constant cocycle, a \Borel compression of $\rho$ over a finite
\Borel subequivalence relation of $E$ gives rise to a \Borel
compression of $\rho$.

\begin{proposition} \label{smooth:compression}
  Suppose that $X$ is a standard \Borel space, $E$ is a
  countable \Borel equivalence relation on $X$, $\rho \from E
  \to \openinterval{0}{\infty}$ is a \Borel cocycle, and there is a
  \Borel compression $\phi \from X \to X$ of $\rho$ over a
  finite \Borel subequivalence relation $F$ of $E$. Then there is
  no $\rho$-invariant \Borel probability measure.
\end{proposition}

\begin{propositionproof}
  By the \Lusin-\Novikov uniformization theorem, there exist a \Borel
  transversal $B \subseteq X$ of $F$, \Borel sets $B_n \subseteq
  B$, and \Borel injections $T_n \from B_n \to X$ for which $\sequence
  {\graph{T_n}}[n \in \N]$ partitions $F \intersection (B \times X)$. If
  $\mu$ is a $\rho$-invariant \Borel measure, then Proposition \ref
  {smooth:injection} ensures that
  \begin{align*}
    \mu(X)
      & = \sum_{n \in \N} \mu(\image{T_n}{B_n}) \\
      & = \sum_{n \in \N} \int_{B_n} \rho(T_n(x), x) \ d\mu(x) \\
      & = \int_B \cardinality[x][\rho]{\equivalenceclass{x}{F}} \ d\mu(x),
  \end{align*}
  whereas Propositions \ref{smooth:injection} and \ref
  {smooth:invariance} imply that
  \begin{align*}
    \mu(X)
      & = \int \cardinality[x][\rho]{\preimage{\phi}{x}} \ d\mu(x) \\
      & = \sum_{n \in \N} \int_{\image{T_n}{B_n}} \cardinality[x][\rho]
        {\preimage{\phi}{x}} \ d\mu(x) \\
      & = \sum_{n \in \N} \int_{B_n} \cardinality[T_n(x)][\rho]{(\inverse
        {\phi} \composition T_n)(x)} \ d(\pushforward{(\inverse{T_n})}
          (\mu))(x) \\
      & = \sum_{n \in \N} \int_{B_n} \cardinality[x][\rho]{(\inverse{\phi}
        \composition T_n)(x)} \ d\mu(x) \\
      & = \int_B \cardinality[x][\rho]{\preimage{\phi}{\equivalenceclass
        {x}{F}}} \ d\mu(x).
  \end{align*}
  As the set $A = \set{x \in B}[{\cardinality[x][\rho]{\preimage{\phi}
  {\equivalenceclass{x}{F}}} < \cardinality[x][\rho]{\equivalenceclass
  {x}{F}}}]$ is $E$-complete, it follows that if $\mu(X) > 0$, then $\mu
  (A) > 0$. As $\cardinality[x][\rho]{\preimage{\phi}{\equivalenceclass
  {x}{F}}} \le \cardinality[x][\rho]{\equivalenceclass{x}{F}}$ for all $x \in
  B$, it follows that if $\mu(A) > 0$, then $\mu(X) = \infty$.
\end{propositionproof}

We next note the useful fact that smoothness is invariant under
quotients by finite \Borel subequivalence relations of $E$.

\begin{proposition} \label{smooth:quotient}
  Suppose that $X$ is a standard \Borel space, $E$ is a countable
  \Borel equivalence relation on $X$, $\rho \from E \to \openinterval
  {0}{\infty}$ is a \Borel cocycle, and $F$ is a finite \Borel
  subequivalence relation of $E$. Then $\rho$ is smooth if and only
  if $\rho / F$ is smooth.
\end{proposition}

\begin{propositionproof}
  By partitioning $X$ into countably-many $F$-invariant \Borel sets,
  we can assume that there is a real number $r > 1$ with
  $\cardinality[x][\rho]{\equivalenceclass{x}{F}} \le r$ for all $x \in X$.
  As $\saturation{Y}{F} / F$ is $\lacunaritygraph{\openinterval{1/r}{r}}
  {\rho / F}$-independent for all $\lacunaritygraph{\openinterval
  {1/r^2}{r^2}}{\rho}$-independent sets $Y \subseteq X$, the
  smoothness of $\rho$ yields that of $\rho / F$. As every
  $F$-invariant set $Y \subseteq X$ for which $Y / F$ is
  $\lacunaritygraph{\openinterval{1/r^2}{r^2}}{\rho / F}$-independent
  is itself $(\lacunaritygraph{\openinterval{1/r}{r}}{\rho} \setminus
  F)$-independent, the smoothness of $\rho / F$ yields that of $\rho$.
\end{propositionproof}

Generalizing the \Dougherty-\Jackson-\Kechris observation that
there is a \Borel compression of $E$ if and only if there is an
aperiodic smooth \Borel subequivalence relation of $E$, we have
the following.

\begin{proposition} \label{smooth:characterization:compression}
  Suppose that $X$ is a standard \Borel space, $E$ is a countable
  \Borel equivalence relation on $X$, and $\rho \from E \to
  \openinterval{0}{\infty}$ is a \Borel cocycle. Then the following are
  equivalent:
  \begin{enumerate}
    \item There is an injective \Borel compression of the quotient of
      $\rho$ by a finite \Borel subequivalence relation of $E$.
    \item There is a \Borel subequivalence relation of $E$ on which
      $\rho$ is aperiodic and smooth.
    \item There exist an $E$-invariant \Borel set $B \subseteq X$ and
      a \Borel subequivalence relation $F$ of $E$ such that
      $\restriction{F}{\setcomplement{B}}$ is smooth, $\restriction{\rho}
      {(\restriction{F}{\setcomplement{B}})}$ is aperiodic, and there is a
      strictly $(\restriction{\rho}{(\restriction{F}{B})})$-increasing \Borel
      automorphism.
  \end{enumerate}
\end{proposition}

\begin{propositionproof}
  To see $(1) \implies (2)$, observe that by Proposition \ref
  {smooth:quotient}, we can assume that there is an injective
  \Borel compression $\phi \from X \to X$ of $\rho$. Set $A = \set
  {x \in X}[{\cardinality[x][\rho]{\preimage{\phi}{x}} < 1}]$, and let
  $F$ be the orbit equivalence relation generated by $\phi$. As the
  sets $A_r = \set{x \in X}[{\cardinality[x][\rho]{\preimage{\phi}{x}} <
  r}]$ are $(\restriction{\rho}{F})$-lacunary for all $r < 1$, it follows
  that $\restriction{\rho}{(\restriction{F}{A})}$ is smooth, thus
  $\restriction{\rho}{(\restriction{F}{\saturation{A}{F}})}$ is aperiodic
  and smooth. By the \Lusin-\Novikov uniformization theorem, there
  is a \Borel extension $\psi \from X \to \saturation{A}{F}$ of the
  identity function on $\saturation{A}{F}$ whose graph is contained in
  $E$, in which case the restriction of $\rho$ to the pullback of
  $\restriction{F}{\saturation{A}{F}}$ through $\psi$ is aperiodic and
  smooth.
  
  To see $(2) \implies (3)$, note that if condition (2) holds, then
  Proposition \ref{smooth:smoothandstrictlyincreasing} immediately
  yields the weakening of condition (3) in which the set $B$ need not
  be $E$-invariant. To see that this weakening yields condition (3)
  itself, note that if $B' \subseteq X$ is a \Borel set and $F'$ is a
  smooth \Borel subequivalence relation of $\restriction{E}{B'}$ for
  which $\restriction{\rho}{F'}$ is aperiodic, then the \Lusin-\Novikov
  uniformization theorem yields a \Borel extension $\pi \from \saturation
  {B'}{E} \to B'$ of the identity function on $B'$ whose graph is
  contained in $E$, the subequivalence relation $F''$ of $\restriction{E}
  {\saturation{B'}{E}}$ given by $x \mathrel{F''} y \iff \pi(x) \mathrel{F'}
  \pi(y)$ is smooth, and $\restriction{\rho}{F''}$ is aperiodic.
  
  It only remains to note that Proposition \ref{smooth:aperiodic}
  yields $(3) \implies (1)$.
\end{propositionproof}

We close this section by noting that, at least when $\rho$ is smooth,
the existence of an injective \Borel compression of the quotient of
$\rho$ by a finite \Borel subequivalence relation of $E$ is the sole
obstacle to the existence of a $\rho$-invariant \Borel probability measure.

\begin{proposition} \label{smooth:characterization}
  Suppose that $X$ is a standard \Borel space, $E$ is a countable
  \Borel equivalence relation on $X$, and $\rho \from E \to
  \openinterval{0}{\infty}$ is a smooth \Borel cocycle. Then exactly
  one of the following holds:
  \begin{enumerate}
    \item There is an injective \Borel compression of the quotient
      of $\rho$ by a finite \Borel subequivalence relation of $E$.
    \item There is a $\rho$-invariant \Borel probability measure.
  \end{enumerate}
\end{proposition}

\begin{propositionproof}
  Proposition \ref{smooth:compression} ensures that conditions
  (1) and (2) are mutually exclusive. To see that at least one of them
  holds, note first that if $\rho$ is aperiodic, then Proposition \ref
  {smooth:aperiodic} yields a finite \Borel subequivalence relation
  $F$ of $E$ for which there is a strictly $(\rho / F)$-increasing
  \Borel injection. And if there is a $\rho$-finite equivalence class $C$
  of $E$, then the \Borel probability measure $\mu$ on $X$, given by
  $\mu(B) = \cardinality[C][\rho]{B \intersection C}$ for all \Borel sets
  $B \subseteq X$, is $\rho$-invariant.
\end{propositionproof}

\section{Coboundaries} \label{coboundaries}

We say that a \Borel cocycle $\rho \from E \to \openinterval{0}{\infty}$
is a \definedterm{\Borel coboundary} if there is a \Borel function $f
\from X \to \openinterval{0}{\infty}$ such that $\rho(x, y) = f(x) / f(y)$
for all $\pair{x}{y} \in E$. The following observation shows that, even
for \Borel coboundaries, the equivalent conditions of Proposition
\ref{smooth:characterization:compression} do not characterize the
non-existence of an invariant \Borel probability measure.

\begin{proposition} \label{coboundaries:example}
  Suppose that $X$ is a standard \Borel space and $E$ is an aperiodic
  countable \Borel equivalence relation on $X$ admitting an
  invariant \Borel probability measure. Then there is a
  \Borel coboundary $\rho \from E \to \openinterval{0}{\infty}$
  with the property that there is neither an injective \Borel
  compression of the quotient of $\rho$ by a finite \Borel
  subequivalence relation of $E$ nor a $\rho$-invariant \Borel
  probability measure.
\end{proposition}

\begin{propositionproof}
  Set $B_0 = X$ and let $\iota_0 \from B_0 \to B_0$ be the identity
  function. Recursively apply \cite[Proposition 7.4]{KechrisMiller}
  to obtain \Borel sets $B_{n+1} \subseteq \image{\iota_n}{B_n}$ and
  \Borel involutions $\iota_{n+1} \from \image{\iota_n}{B_n} \to \image
  {\iota_n}{B_n}$ such that the graph of $\iota_{n+1}$ is contained in
  $E$ and the sets $B_{n+1}$ and $\image{\iota_{n+1}}{B_{n+1}}$
  partition $\image{\iota_n}{B_n}$ for all $n \in \N$. For each $x \in X$,
  let $n(x)$ be the maximal natural number for which $x \in B_{n(x)}$,
  and set $f(x) = 2^{n(x)}$. Define $\rho \from E \to \openinterval{0}
  {\infty}$ by setting $\rho(x, y) = f(x) / f(y)$ for all $\pair{x}{y} \in E$.
  
  To see that there is no $\rho$-invariant \Borel probability measure,
  note that if $\mu$ is a $\rho$-invariant \Borel measure, then the
  fact that $\image{\iota_{n+1}}{B_{n+2}}$ and $\image{(\iota_{n+1}
  \composition \iota_{n+2})}{B_{n+2}}$ partition $B_{n+1}$ for all $n
  \in \N$ ensures that
  \begin{equation*}
    \mu(B_{n+1})
      = \int_{B_{n+2}} \rho(\iota_{n+1}(x), x) + \rho((\iota_{n+1}
        \composition \iota_{n+2})(x), x) \ d\mu(x)
      = \mu(B_{n+2})
  \end{equation*}
  for all $n \in \N$, thus $\mu(X) \in \set{0, \infty}$.
  
  Suppose, towards a contradiction, that there is an injective
  \Borel compression of the quotient of $\rho$ by a finite \Borel
  subequivalence relation of $E$. Then Proposition \ref
  {smooth:characterization:compression} yields an $E$-invariant
  \Borel set $A \subseteq X$ and a \Borel subequivalence relation
  $F$ of $E$ such that $\restriction{F}{\setcomplement{A}}$ is
  smooth, $\restriction{\rho}{(\restriction{F}{\setcomplement{A}})}$
  is aperiodic, and there is a strictly $(\restriction{\rho}
  {(\restriction{F}{A})})$-increasing \Borel automorphism $\phi
  \from A \to A$. Fix an $E$-invariant \Borel probability measure
  $\mu$. As $\image{\iota_n}{B_{n+1}}$ and $\image{(\iota_n
  \composition \iota_{n+1})}{B_{n+1}}$ partition $B_n$ for all $n \in
  \N$, it follows that $\mu(B_n) = 2 \mu(B_{n+1})$ for all $n \in \N$.
  As the aperiodicity of $\restriction{\rho}{(\restriction{F}
  {\setcomplement{A}})}$ yields that of $\restriction{F}
  {\setcomplement{A}}$, Propositions \ref{smooth:aperiodic} and
  \ref{smooth:compression} imply that $A$ is $\mu$-conull, thus so
  too is $A \intersection \union[n \in \N][B_{n+1}]$. As the definition
  of $\rho$ ensures that $\image{\phi}{A \intersection \union[n \in \N]
  [B_{n+1}]} \subseteq A \intersection \union[n \in \N][B_{n+2}]$, and
  the latter set has $\mu$-measure $1/2$, this contradicts
  $E$-invariance.
\end{propositionproof}

The following fact yields an equivalent of $\rho$-invariance that
will prove useful when dealing with finite \Borel subequivalence
relations.

\begin{proposition} \label{coboundaries:integral}
  Suppose that $X$ is a standard \Borel space, $E$ is a countable
  \Borel equivalence relation on $X$, $\rho \from E \to \openinterval
  {0}{\infty}$ is a \Borel cocycle, and $\mu$ is a $\rho$-invariant
  \Borel measure on $X$. Then $\mu(B) = \int \cardinality
  [\equivalenceclass{x}{F}][\rho]{B \intersection \equivalenceclass{x}
  {F}} \ d\mu(x)$ for all \Borel sets $B \subseteq X$ and finite \Borel
  subequivalence relations $F$ of $E$.
\end{proposition}

\begin{propositionproof}
  Fix a \Borel transversal $A \subseteq X$ of $F$, \Borel sets $A_n
  \subseteq A$, and \Borel injections $T_n \from A_n \to X$ with the
  property that $\sequence{\graph{T_n}}[n \in \N]$ partitions $F
  \intersection (A \times X)$, and observe that
  \begin{align*}
    \int \cardinality[\equivalenceclass{x}{F}][\rho]{B \intersection
      \equivalenceclass{x}{F}} \ d\mu(x)
      & = \sum_{n \in \N} \int_{\image{T_n}{A_n}} \cardinality
        [\equivalenceclass{x}{F}][\rho]{B \intersection
          \equivalenceclass{x}{F}} \ d\mu(x) \\
      & = \sum_{n \in \N} \int_{A_n} \cardinality[\equivalenceclass
        {x}{F}][\rho]{B \intersection \equivalenceclass{x}{F}}
          \ d(\pushforward{(\inverse{T_n})}{\mu})(x) \\
      & = \sum_{n \in \N} \int_{A_n} \cardinality
        [\equivalenceclass{x}{F}][\rho]{B \intersection 
          \equivalenceclass{x}{F}} \rho(T_n(x), x) \ d\mu(x) \\
      & = \int_A \cardinality[x][\rho]{B \intersection 
        \equivalenceclass{x}{F}} \ d\mu(x) \\
      & = \sum_{n \in \N} \int_{A_n \intersection \preimage{T_n}{B}}
        \rho(T_n(x), x) \ d\mu(x) \\
      & = \sum_{n \in \N} \mu(\image{T_n}{A_n} \intersection B) \\
      & = \mu(B)
  \end{align*}
  by Proposition \ref{smooth:injection}.
\end{propositionproof}

Given a \Borel set $R \subseteq X \times X$ with countable vertical
sections and a \Borel function $\rho \from R \to \openinterval{0}
{\infty}$, we say that a \Borel measure $\mu$ on $X$ is \definedterm
{$\rho$-invariant} if $\mu(\image{T}{B}) = \int_B \rho(T(x), x)
\ d\mu(x)$ for all \Borel sets $B \subseteq X$ and \Borel injections
$T \from B \to X$ whose graphs are contained in $\inverse{R}$. The
\definedterm{composition} of sets $R \subseteq X \times Y$ and $S
\subseteq Y \times Z$ is given by $R \composition S = \set{\pair{x}
{z} \in X \times Z}[\exists y \in Y \ x \mathrel{R} y \mathrel{S} z]$.
The \Lusin-\Novikov uniformization theorem ensures that if $R$ and
$S$ are \Borel sets with countable vertical sections, then so too is
their composition. The following fact will prove useful in verifying
$\rho$-invariance.

\begin{proposition} \label{coboundaries:composition}
  Suppose that $X$ is a standard \Borel space, $E$ is a countable
  \Borel equivalence relation on $X$, $R, S \subseteq E$ are \Borel,
  and $\rho \from E \to \openinterval{0}{\infty}$ is a \Borel cocycle.
  Then every $(\restriction{\rho}{(R \union S)})$-invariant \Borel
  measure $\mu$ is $(\restriction{\rho}{(R \composition
  S)})$-invariant.
\end{proposition}

\begin{propositionproof}
  Note first that if $B \subseteq X$ is a \Borel set, $T_S \from B \to
  X$ is a \Borel injection whose graph is contained in $\inverse{S}$,
  and $T_R \from \image{T_S}{B} \to X$ is a \Borel injection whose
  graph is contained in $\inverse{R}$, then
  \begin{align*}
    \mu(\image{(T_R \composition T_S)}{B})
      & = \int_{\image{T_S}{B}} \rho(T_R(x), x) \ d\mu(x) \\
      & = \int_B \rho((T_R \composition T_S)(x), T_S(x))
        \ d(\pushforward{(\inverse{T_S})}
          {\mu})(x) \\
      & = \int_B \rho((T_R \composition T_S)(x), x) \ d\mu(x).
  \end{align*}
  As the \Lusin-\Novikov uniformization theorem ensures that every
  \Borel injection whose graph is contained in $\inverse{(R
  \composition S)}$ can be decomposed into a disjoint union of
  countably-many \Borel injections of the form $T_R \composition
  T_S$ as above, the proposition follows.
\end{propositionproof}

We say that \Borel cocycles $\rho \from E \to \openinterval{0}{\infty}$
and $\sigma \from E \to \openinterval{0}{\infty}$ are \definedterm
{\Borel cohomologous} if their ratio is a \Borel coboundary. We say
that a \Borel function $f \from X \to \openinterval{0}{\infty}$
\definedterm{witnesses} that $\rho$ and $\sigma$ are \Borel
cohomologous if $f(x) / f(y) = \sigma(x, y) / \rho(x, y)$ for all $\pair
{x}{y} \in E$.

\begin{proposition} \label{coboundaries:invariantconversion}
  Suppose that $X$ is a standard \Borel space, $E$ is a countable
  \Borel equivalence relation on $X$, $f \from X \to \openinterval
  {0}{\infty}$ is a \Borel function witnessing that \Borel cocycles
  $\rho, \sigma \from E \to \openinterval{0}{\infty}$ are \Borel
  cohomologous, and $\mu$ is a $\rho$-invariant \Borel measure.
  Then the \Borel measure given by $\nu(B) = \int_B f \ d\mu$ is
  $\sigma$-invariant.
\end{proposition}

\begin{propositionproof}
  Simply observe that if $B \subseteq X$ is a \Borel set and $T \from
  X \to X$ is a \Borel automorphism whose graph is contained in $E$,
  then
  \begin{align*}
    \nu(\image{T}{B})
      & = \int_{\image{T}{B}} f \ d\mu \\
      & = \int_B f \composition T \ d(\pushforward{(\inverse{T})}{\mu}) \\
      & = \int_B (f \composition T)(x) \rho(T(x), x) \ d\mu(x) \\
      & = \int_B f(x) \sigma(T(x), x) \ d\mu(x) \\
      & = \int_B \sigma(T(x), x) \ d\nu(x)
  \end{align*}
  by $\rho$-invariance.
\end{propositionproof}

We say that a \Borel set $B \subseteq X$ has \definedterm
{$\rho$-density} at least $\epsilon$ if there is a finite \Borel
subequivalence relation $F$ of $E$ such that $\cardinality
[\equivalenceclass{x}{F}][\rho]{B \intersection \equivalenceclass{x}
{F}} \ge \epsilon$ for all $x \in X$. We say that a \Borel set $B
\subseteq X$ has \definedterm{positive $\rho$-density} if there
exists $\epsilon > 0$ for which $B$ has $\rho$-density at least
$\epsilon$.

\begin{proposition} \label{coboundaries:extension:measure}
  Suppose that $X$ is a standard \Borel space, $E$ is a countable
  \Borel equivalence relation on $X$, $\rho \from E \to \openinterval
  {0}{\infty}$ is a \Borel cocycle, and $B \subseteq X$ is a \Borel
  set with positive $\rho$-density. Then every $(\restriction{\rho}
  {(\restriction{E}{B})})$-invariant finite \Borel measure $\mu$
  extends to a $\rho$-invariant finite \Borel measure.
\end{proposition}

\begin{propositionproof}
  Fix $\epsilon > 0$ for which $B$ has $\rho$-density at least
  $\epsilon$, as well as a finite \Borel subequivalence relation $F$ of
  $E$ such that $\cardinality[\equivalenceclass{x}{F}][\rho]{B
  \intersection \equivalenceclass{x}{F}} \ge \epsilon$ for all $x \in
  X$, and let $\extension{\mu}$ be the \Borel measure on $X$ given by
  \begin{equation*}
    \extension{\mu}(A) = \int \cardinality[B \intersection
      \equivalenceclass{x}{F}][\rho]{A \intersection \equivalenceclass{x}
        {F}} \ d\mu(x)
  \end{equation*}
  for all \Borel sets $A \subseteq X$.
  
  As $\extension{\mu}(X) \le \mu(B) / \epsilon$, it follows that
  $\extension{\mu}$ is finite, and Proposition \ref
  {coboundaries:integral} ensures that $\mu = \restriction{\extension
  {\mu}}{B}$.
  
  \begin{lemma} \label{coboundaries:extension:measure:integral}
    Suppose that $f \from X \to \closedopeninterval{0}{\infty}$ is a
    \Borel function. Then $\int f \ d\extension{\mu} = \int \sum_{y \in
    \equivalenceclass{x}{F}} f(y) \cardinality[B \intersection
    \equivalenceclass{x}{F}][\rho]{\set{y}} \ d\mu(x)$.
  \end{lemma}
  
  \begin{lemmaproof}
    It is sufficient to check the special case that $f$ is the characteristic
    function of a \Borel set, which is a direct consequence of the
    definition of $\extension{\mu}$.
  \end{lemmaproof}
  
  \begin{lemma}
    The measure $\extension{\mu}$ is $(\restriction{\rho}{F})$-invariant.
  \end{lemma}
  
  \begin{lemmaproof}
    Simply observe that if $A \subseteq X$ is a \Borel set and
    $T \from X \to X$ is a \Borel automorphism whose graph is
    contained in $F$, then
    \begin{align*}
      \int_A \rho(T(x), x) \ d\extension{\mu}(x)
        & = \int \sum_{y \in A \intersection \equivalenceclass{x}{F}}
          \rho(T(y), y) \cardinality[B \intersection \equivalenceclass{x}
            {F}][\rho]{\set{y}} \ d\mu(x) \\
        & = \int \cardinality[B \intersection \equivalenceclass{x}{F}]
          [\rho]{\image{T}{A \intersection \equivalenceclass{x}{F}}}
            \ d\mu(x) \\
        & = \extension{\mu}(\image{T}{A})
    \end{align*}
    by Lemma \ref{coboundaries:extension:measure:integral}.
  \end{lemmaproof}
  
  As $E = F \composition (E \intersection (B \times B)) \composition
  F$, two applications of Proposition \ref{coboundaries:composition}
  ensure that $\extension{\mu}$ is $\rho$-invariant.
\end{propositionproof}

The primary argument of this section will hinge on the following
approximation lemma.

\begin{proposition} \label{coboundaries:approximation}
  Suppose that $X$ is a standard \Borel space, $E$ is a countable
  \Borel equivalence relation on $X$, and $\rho \from E \to
  \openinterval{0}{\infty}$ is a \Borel cocycle. Then for all \Borel sets
  $A \subseteq X$ and positive real numbers $r < 1$, there exist an
  $E$-invariant \Borel set $B \subseteq X$, a \Borel set $C \subseteq
  B$, and a finite \Borel subequivalence relation $F$ of $\restriction
  {E}{C}$ such that $\restriction{\rho}{(\restriction{E}{\setcomplement
  {B}})}$ is smooth, $r < \cardinality[\equivalenceclass{x}{F} \setminus
  A][\rho]{A \intersection \equivalenceclass{x}{F}} < 1$ for all $x \in
  C$, and $A \intersection \equivalenceclass{x}{E} \subseteq C$ or
  $\equivalenceclass{x}{E} \setminus A \subseteq C$ for all $x \in B$.
\end{proposition}

\begin{propositionproof}
  By \cite[Lemma 7.3]{KechrisMiller}, there is a maximal \Borel set
  $\calS$ of pairwise disjoint non-empty finite sets $S \subseteq X$
  with $S \times S \subseteq E$ and $r < \cardinality[S
  \setminus A][\rho]{A \intersection S} < 1$. Set $D =  A \setminus
  \union[\calS]$ and $D' = (\setcomplement{A}) \setminus \union
  [\calS]$.

  \begin{lemma} \label{coboundaries:approximation:bound}
    Suppose that $\pair{x}{x'} \in E$. Then there exists a real number
    $s > 1$ with the property that $x$ has only finitely-many
    $\lacunaritygraph{\openinterval{1/s}{s}}{\rho}$-neighbors in $D$
    or $x'$ has only finitely-many $\lacunaritygraph{\openinterval{1/s}
    {s}}{\rho}$-neighbors in $D'$.
  \end{lemma}
  
  \begin{lemmaproof}
    Fix $n, n' \in \N$ such that $(n / n') \rho(x, x')$ lies strictly between
    $r$ and $1$, and fix $s > 1$ sufficiently small that $(n / n') \rho(x,
    x')$ lies strictly between $r s^2$ and $1 / s^2$. Suppose, towards
    a contradiction, that there are sets $S \subseteq D$ and $S'
    \subseteq D'$ of $\lacunaritygraph{\openinterval{1/s}{s}}
    {\rho}$-neighbors of $x$ and $x'$ of cardinalities $n$ and $n'$.
    Then $n / s < \cardinality[x][\rho]{S} < n s$ and $n' \rho(x', x) / s <
    \cardinality[x][\rho]{S'} < n' \rho(x', x) s$, so the $\rho$-size of $S$
    relative to $S'$ lies strictly between $(n / n') \rho(x, x') / s^2$ and
    $(n / n') \rho(x, x') s^2$. As these bounds lie strictly between $r$
    and $1$, this contradicts the maximality of $\calS$.
  \end{lemmaproof}
  
  Lemma \ref{coboundaries:approximation:bound} ensures that
  $\saturation{D}{E} \intersection \saturation{D'}{E}$ is contained in the
  $E$-saturation of the union of the sets of the form $\set{x \in D}
  [\cardinality{D \intersection \verticalsection{(\lacunaritygraph
  {\openinterval{1/s}{s}}{\rho})}{x}} < \aleph_0]$ and $\set{x \in D'}
  [\cardinality{D' \intersection \verticalsection{(\lacunaritygraph
  {\openinterval{1/s}{s}}{\rho})}{x}} < \aleph_0]$, so $\restriction{\rho}
  {(\restriction{E}{(\saturation{D}{E} \intersection \saturation{D'}{E})})}$
  is smooth. Set $B = \setcomplement{(\saturation{D}{E} \intersection
  \saturation{D'}{E})}$ and $C = B \intersection \union[\calS]$, and let
  $F$ be the equivalence relation on $C$ whose classes are the
  subsets of $C$ in $\calS$.
\end{propositionproof}

We say that a \Borel set $B \subseteq X$ has \definedterm
{$\sigma$-positive $\rho$-density} if $X$ is the union of
countably-many $E$-invariant \Borel sets $A_n \subseteq X$ for
which $A_n \intersection B$ has positive $(\restriction{\rho}
{(\restriction{E}{A_n})})$-density.

\begin{theorem} \label{coboundaries:extension}
  Suppose that $X$ is a standard \Borel space, $E$ is a countable
  \Borel equivalence relation on $X$, $\rho \from E \to \openinterval
  {0}{\infty}$ is a \Borel cocycle, and $A \subseteq X$ is an
  $E$-complete \Borel set. Then $X$ is the union of an $E$-invariant
  \Borel set $B \subseteq X$ for which $\restriction{\rho}{(\restriction
  {E}{B})}$ is smooth, an $E$-invariant \Borel set $C \subseteq X$
  for which $A \intersection C$ has $\sigma$-positive $(\restriction
  {\rho}{(\restriction{E}{C})})$-density, and an $E$-invariant \Borel set
  $D \subseteq X$ for which there is a finite-to-one \Borel
  compression of the quotient of $\restriction{\rho}{(\restriction{E}
  {D})}$ by a finite \Borel subequivalence relation of $\restriction{E}
  {D}$.
\end{theorem}

\begin{theoremproof}
  Fix a positive real number $r < 1$. We will show that, after throwing
  out countably-many $E$-invariant \Borel sets $B \subseteq X$ for
  which $\restriction{\rho}{(\restriction{E}{B})}$ is smooth, as well as
  countably-many $E$-invariant \Borel sets $C \subseteq X$ for
  which $A \intersection C$ has positive $(\restriction{\rho}{(\restriction
  {E}{C})})$-density, there are increasing sequences of finite \Borel
  subequivalence relations $F_n$ of $E$ and $E$-complete
  $F_n$-invariant \Borel sets $A_n \subseteq X$ with the property
  that $r < \cardinality[(A_{n+1} \setminus A_n) \intersection
  \equivalenceclass{x}{F_{n+1}}][\rho]{A_n \intersection
  \equivalenceclass{x}{F_{n+1}}} < 1$ for all $n \in \N$ and $x \in
  A_n$.
  
  We begin by setting $A_0 = A$ and letting $F_0$ be equality.
  Suppose now that $n \in \N$ and we have already found $A_n$
  and $F_n$. By applying Proposition \ref
  {coboundaries:approximation} to $A_n / F_n$, and throwing out an
  $E$-invariant \Borel set $B \subseteq X$ for which $\restriction
  {\rho}{(\restriction{E}{B})}$ is smooth, we obtain a finite \Borel
  subequivalence relation $F_{n+1} \supseteq F_n$ of $E$ and an
  $F_{n+1}$-invariant \Borel set $A_{n+1} \subseteq X$ such that $r
  < \cardinality[\equivalenceclass{x}{F_{n+1}} \setminus A_n][\rho]
  {A_n \intersection \equivalenceclass{x}{F_{n + 1}}} < 1$  for all $x
  \in A_{n+1}$, and $A_n \intersection \equivalenceclass{x}{E}
  \subseteq A_{n+1}$ or $\equivalenceclass{x}{E} \setminus A_n
  \subseteq A_{n+1}$ for all $x \in X$. By throwing out an
  $E$-invariant \Borel set $C \subseteq X$ for which $A \intersection
  C$ has positive $(\restriction{\rho}{(\restriction{E}{C})})$-density, we
  can assume that $A_n \subseteq A_{n+1}$, completing the recursive
  construction.
  
  Set $B_n = A_n \setminus \union[m < n][A_m]$ and define $\phi_n
  \from B_n / F_n \to B_{n+1} / F_{n+1}$ by setting $\phi_n(B_n
  \intersection \equivalenceclass{x}{F_n}) = B_{n+1} \intersection
  \equivalenceclass{x}{F_{n+1}}$ for all $n \in \N$ and $x \in B_n$.
  Then the union of $\union[n \in \N][\phi_n]$ and the identity function
  on $\setcomplement{\union[n \in \N][A_n]}$ is a \Borel compression
  of the quotient of $\rho$ by the union of $\union[n \in \N][\restriction
  {F_n}{B_n}]$ and equality.
\end{theoremproof}

As a corollary, we obtain the desired characterization.

\begin{theorem} \label{coboundaries:compression}
  Suppose that $X$ is a standard \Borel space, $E$ is a countable
  \Borel equivalence relation on $X$, and $\rho \from E \to
  \openinterval{0}{\infty}$ is a \Borel coboundary. Then exactly one
  of the following holds:
  \begin{enumerate}
    \item There is a finite-to-one \Borel compression of the quotient of
      $\rho$ by a finite \Borel subequivalence relation of $E$.
    \item There is a $\rho$-invariant \Borel probability measure.
  \end{enumerate}
\end{theorem}

\begin{theoremproof}
  Proposition \ref{smooth:compression} ensures that conditions (1)
  and (2) are mutually exclusive. To see that at least one of them
  holds, fix a bounded open neighborhood $U \subseteq
  \openinterval{0}{\infty}$ of $1$. As $\rho$ is a \Borel coboundary,
  the \Lusin-\Novikov uniformization theorem implies that there
  is an $E$-complete \Borel set $A \subseteq X$ for which $\image
  {\rho}{\restriction{E}{A}} \subseteq U$. By Theorem \ref
  {coboundaries:extension}, after throwing out $E$-invariant \Borel
  sets $B \subseteq X$ and $D \subseteq X$ for which $\restriction
  {\rho}{(\restriction{E}{B})}$ is smooth and there is a finite-to-one
  \Borel compression of the quotient of $\restriction{\rho}{(\restriction
  {E}{D})}$ by a finite \Borel subequivalence relation of $\restriction
  {E}{D}$, we can assume that $A$ has $\sigma$-positive
  $\rho$-density.
  
  If there is a $(\restriction{\rho}{(\restriction{E}{A})})$-invariant \Borel
  probability measure $\mu$, then by passing to an $(\restriction{E}
  {A})$-invariant $\mu$-positive \Borel set, we can assume that
  $A$ has positive $\rho$-density, in which case Proposition \ref
  {coboundaries:extension:measure} yields a $\rho$-invariant
  \Borel probability measure.
  
  If there is no $(\restriction{\rho}{(\restriction{E}{A})})$-invariant
  \Borel probability measure, then Proposition \ref
  {coboundaries:invariantconversion} ensures that there is no
  $(\restriction{E}{A})$-invariant \Borel probability measure, in
  which case the \Becker-\Kechris generalization of \Nadkarni's
  theorem and the \Dougherty-\Jackson-\Kechris characterization
  of the existence of a \Borel compression yield an aperiodic smooth
  \Borel subequivalence relation $F$ of $\restriction{E}{A}$. Then
  $\restriction{\rho}{F}$ is smooth, and the fact that $\restriction{\rho}
  {(\restriction{E}{A})}$ is bounded ensures that $\restriction{\rho}{F}$
  is also aperiodic. Fix a \Borel extension $\phi \from X \to A$ of the
  identity function on $A$ whose graph is contained in $E$, and
  observe that $\rho$ is aperiodic and smooth on the pullback of $F$
  through $\phi$. Proposition \ref{smooth:aperiodic} therefore yields
  an injective \Borel compression of the quotient of $\rho$ by a finite
  \Borel subequivalence relation of $E$.
\end{theoremproof}

\section{The general case} \label{cocycles}

Here we generalize \Nadkarni's theorem to \Borel cocycles. As in
\S\ref{coboundaries}, our primary argument will hinge on a pair of
approximation lemmas. Given a finite set $S \subseteq X$ for which
$S \times S \subseteq E$, let $\measure{S}{\rho}$ be the \Borel
probability measure on $X$ given by $\measure{S}{\rho}(B) =
\cardinality[S][\rho]{B \intersection S}$.

\begin{proposition} \label{cocycles:approximation}
  Suppose that $X$ is a standard \Borel space, $E$ is a countable
  \Borel equivalence relation on $X$, $\rho \from E \to \openinterval
  {0}{\infty}$ is a \Borel cocycle, $f \from X \to \closedopeninterval{0}
  {\infty}$ is \Borel, $\delta > 0$, and $\epsilon > \sup_{\pair{x}{y} \in
  E} f(x) - f(y)$. Then there exist an $E$-invariant \Borel set $B
  \subseteq X$ and a finite \Borel subequivalence relation $F$ of
  $\restriction{E}{B}$ for which $\restriction{\rho}{(\restriction{E}
  {\setcomplement{B}})}$ is smooth and $\delta \epsilon > \sup_{\pair
  {x}{y} \in \restriction{E}{B}} \int f \ d\measure{\equivalenceclass{x}{F}}
  {\rho} - \int f \ d\measure{\equivalenceclass{y}{F}}{\rho}$.
\end{proposition}

\begin{propositionproof}
  We can clearly assume that $\delta < 1$, and since one can
  repeatedly apply the corresponding special case of the proposition
  over the corresponding quotients, we can also assume that $\delta >
  2 / 3$. For each $x \in X$, let $\average{f}(\equivalenceclass{x}{E})$
  be the average of $\inf \image{f}{\equivalenceclass{x}{E}}$ and $\sup
  \image{f}{\equivalenceclass{x}{E}}$. By \cite[Lemma 7.3]
  {KechrisMiller}, there is a maximal \Borel set $\calS$ of pairwise
  disjoint non-empty finite sets $S \subseteq X$ with $S \times S
  \subseteq E$ and $\epsilon(\delta - 1/2) > \absolutevalue{\int f
  \ d\measure{S}{\rho} - \average{f}(\saturation{S}{E})}$. Set $C = \set
  {x \in \setcomplement{\union[\calS]}}[f(x) < \average{f}
  (\equivalenceclass{x}{E})]$ and $D = \set{x \in \setcomplement{\union
  [\calS]}}[f(x) > \average{f}(\equivalenceclass{x}{E})]$.

  \begin{lemma} \label{cocycles:approximation:bound}
    Suppose that $\pair{x}{y} \in E$. Then there exists a real number
    $r > 1$ such that $x$ has only finitely-many $\lacunaritygraph
    {\openinterval{1/r}{r}}{\rho}$-neighbors in $C$ or $y$ has only
    finitely-many $\lacunaritygraph{\openinterval{1/r}{r}}
    {\rho}$-neighbors in $D$.
  \end{lemma}
  
  \begin{lemmaproof}
    As $\delta > 2/3$, a trivial calculation reveals that $-\epsilon (\delta - 1/2)$
    is strictly below the average of $-\epsilon / 2$ and $\epsilon (\delta -
    1/2)$, and that the average of $-\epsilon(\delta - 1/2)$ and $\epsilon /
    2$ is strictly below $\epsilon(\delta - 1/2)$. In particular, by choosing
    $m, n \in \N$ for which the ratios $s = m / (m + n \rho(y, x))$ and
    $t = n \rho(y, x) / (m + n \rho(y, x))$ are sufficiently close to $1/2$,
    we can therefore ensure that the sums $s(\average{f}
    (\equivalenceclass{x}{E}) - \epsilon / 2) + t(\average{f}
    (\equivalenceclass{x}{E}) + \epsilon(\delta - 1/2))$ and $s(\average
    {f}(\equivalenceclass{x}{E}) - \epsilon(\delta - 1/2)) + t(\average
    {f}(\equivalenceclass{x}{E}) + \epsilon / 2)$ both lie strictly
    between $\average{f}(\equivalenceclass{x}{E}) - \epsilon(\delta -
    1/2)$ and $\average{f}(\equivalenceclass{x}{E}) + \epsilon(\delta -
    1/2)$. Fix $r > 1$ such that they lie strictly between
    $(\average{f}(\equivalenceclass{x}{E}) - \epsilon(\delta -
    1/2)) r^2$ and $(\average{f}(\equivalenceclass{x}{E}) + \epsilon(\delta -
    1/2)) / r^2$.
  
    Suppose, towards a contradiction, that there exist sets $S
    \subseteq C$ and $T \subseteq D$ of $\lacunaritygraph
    {\openinterval{1/r}{r}}{\rho}$-neighbors of $x$ and $y$ of
    cardinalities $m$ and $n$. Then $m / r < \cardinality[x][\rho]{S}
    < mr$ and $n \rho(y, x) / r < \cardinality[x][\rho]{T} < n \rho(y,
    x)r$, from which a trivial calculation reveals that $s / r^2 <
    \cardinality[x][\rho]{S} / \cardinality[x][\rho]{S \union T} <
    s r^2$ and $ t / r^2 < \cardinality[x][\rho]{T} / \cardinality[x]
    [\rho]{S \union T} < t r^2$. As $\int f \ d\measure{S}{\rho}$
    lies between $\average{f}(\equivalenceclass{x}{E}) -
    \epsilon / 2$ and $\average{f}(\equivalenceclass{x}{E}) - \epsilon
    (\delta - 1/2)$, and $\int f \ d\measure{T}{\rho}$ lies between
    $\average{f}(\equivalenceclass{x}{E}) + \epsilon(\delta - 1/2)$ and
    $\average{f}(\equivalenceclass{x}{E}) + \epsilon / 2$, it follows that
    $\int f \ d\measure{S \union T}{\rho}$ lies between $(s
    (\average{f}(\equivalenceclass{x}{E}) - \epsilon / 2) + t
    (\average{f}(\equivalenceclass{x}{E}) + \epsilon(\delta - 1/2))) / r^2$
    and $(s(\average{f}(\equivalenceclass{x}{E}) - \epsilon(\delta - 1/2))
    + t (\average{f}(\equivalenceclass{x}{E}) + \epsilon / 2)) r^2$, so
    strictly between $\average{f}(\equivalenceclass{x}{E}) - \epsilon
    (\delta - 1/2)$ and $\average{f}(\equivalenceclass{x}{E}) + \epsilon
    (\delta - 1/2)$, contradicting the maximality of $\calS$.
  \end{lemmaproof}
  
  Lemma \ref{cocycles:approximation:bound} ensures that $\saturation
  {C}{E} \intersection \saturation{D}{E}$ is contained in the
  $E$-saturation of the union of the sets of the form $\set{x \in C}
  [\cardinality{C \intersection \verticalsection{(\lacunaritygraph
  {\openinterval{1/r}{r}}{\rho})}{x}} < \aleph_0]$ and $\set{x \in D}
  [\cardinality{D \intersection \verticalsection{(\lacunaritygraph
  {\openinterval{1/r}{r}}{\rho})}{x}} < \aleph_0]$, so $\restriction{\rho}
  {(\restriction{E}{(\saturation{C}{E} \intersection \saturation{D}{E})})}$
  is smooth. Set $B = \setcomplement{(\saturation{C}{E} \intersection
  \saturation{D}{E})}$, and let $F$ be the equivalence relation on $B$
  whose classes are the subsets of $B$ in $\calS$ together with the
  singletons contained in $B \setminus \union[\calS]$.
\end{propositionproof}

\begin{proposition} \label{cocycles:injection}
  Suppose that $X$ is a standard \Borel space, $E$ is a countable
  \Borel equivalence relation on $X$, $\rho \from E \to \openinterval
  {0}{\infty}$ is a \Borel cocycle, $f, g \from X \to \closedopeninterval
  {0}{\infty}$ are \Borel, and $r > 1$. Then there exist an
  $E$-invariant \Borel set $B \subseteq X$, a \Borel set $C \subseteq
  B$, and a finite \Borel subequivalence relation $F$ of $\restriction
  {E}{B}$ such that $\restriction{\rho}{(\restriction{E}{\setcomplement
  {B}})}$ is smooth and $\int_C f \ d\measure{\equivalenceclass{x}{F}}
  {\rho} \le \int_{B \setminus C} g \ d\measure{\equivalenceclass{x}
  {F}}{\rho} \le r \int_C f \ d\measure{\equivalenceclass{x}{F}}{\rho}$
  for all $x \in B$.
\end{proposition}

\begin{propositionproof}
  As the proposition holds trivially on $\preimage{f}{0} \union
  \preimage{g}{0}$, we can assume that $f, g \from X \to \openinterval
  {0}{\infty}$. By \cite[Lemma 7.3]{KechrisMiller}, there is a maximal
  \Borel set $\calS$ of pairwise disjoint non-empty finite sets $S
  \subseteq X$ with $S \times S \subseteq E$ and $1 < \int_{S
  \setminus T} g \ d\measure{S}{\rho} \mathrel{/} \int_T f \ d\measure
  {S}{\rho} < r$ for some $T \subseteq S$.
  
  Set $D_{U, V} = (\preimage{f}{U} \intersection \preimage{g}{V})
  \setminus \union[\calS]$ for all $U, V \subseteq \openinterval{0}
  {\infty}$.
  
  \begin{lemma} \label{cocycles:injection:bound}
    For all $x \in X$, there exists $s > 1$ such that $x$ has only
    finitely-many $\lacunaritygraph{\openinterval{1/s}{s}}
    {\rho}$-neighbors in $D_{\openinterval{f(x)/s}{f(x)s}, \openinterval
    {g(x)/s}{g(x)s}}$.
  \end{lemma}
  
  \begin{lemmaproof}
    Fix $m, n \in \N$ for which $1 < (g(x) / f(x))(n / m) < r$, as well as
    $s > 1$ sufficiently large that $s^6 < (g(x) / f(x))(n / m) < r / s^6$.
    Suppose, towards a contradiction, that there is a set $S \subseteq
    D_{\openinterval{f(x)/s}{f(x)s}, \openinterval{g(x)/s}{g(x)s}}$ of
    $\lacunaritygraph{\openinterval{1/s}{s}}{\rho}$-neighbors of $x$ of
    cardinality $k = m + n$, and fix $T \subseteq S$ of cardinality
    $m$. Then $f(x) \measure{S}{\rho}(T) / s < \int_T f \ d\measure{S}
    {\rho} < f(x) \measure{S}{\rho}(T) s$ and $(m / k) / s^2 <
    \measure{S}{\rho}(T) < (m / k) s^2$, so $f(x) (m / k) / s^3 < \int_T
    f \ d\measure{S}{\rho} < f(x) (m / k) s^3$. And $g(x) \measure{S}
    {\rho}(S \setminus T) / s < \int_{S \setminus T} g \ d\measure{S}
    {\rho} < g(x) \measure{S}{\rho}(S \setminus T) s$ and $(n / k) /
    s^2 < \measure{S}{\rho}(S \setminus T) < (n / k) s^2$, so $g(x) (n
    / k) / s^3 < \int_{S \setminus T} g \ d\measure{S}{\rho} < g(x) (n /
    k) s^3$. It follows that $\int_{S \setminus T} g \ d\measure{S}{\rho}
    \mathrel{/} \int_T f \ d\measure{S}{\rho}$ lies strictly between
    $(g(x) / f(x))(n / m) / s^6$ and $(g(x) / f(x))(n / m) s^6$, and
    therefore strictly between $1$ and $r$, contradicting the
    maximality of $\calS$.
  \end{lemmaproof}
  
  As Lemma \ref{cocycles:approximation:bound} ensures that
  $\setcomplement{\union[\calS]}$ is contained in the union of the sets
  of the form $\set{x \in D_{U,V}}[\cardinality{D_{U,V} \intersection
  \verticalsection{(\lacunaritygraph{\openinterval{1/s}{s}}{\rho})}{x}} <
  \aleph_0]$, it follows that $\restriction{\rho}{(\restriction{E}{\saturation
  {\setcomplement{\union[\calS]}}{E}})}$ is smooth. Set $B = 
  \setcomplement{\saturation{\setcomplement{\union[\calS]}}{E}}$, let
  $F$ be the \Borel equivalence relation on $B$ whose classes are the
  subsets of $B$ in $\calS$, and appeal to the \Lusin-\Novikov
  uniformization theorem to obtain a \Borel set $C \subseteq B$ with
  the property that $1 < \int_{B \setminus C} g \ d\measure
  {\equivalenceclass{x}{F}}{\rho} \mathrel{/} \int_C f \ d\measure
  {\equivalenceclass{x}{F}}{\rho}< r$ for all $x \in B$.
\end{propositionproof}

We are now ready to establish our primary result.

\begin{theorem}
  Suppose that $X$ is a standard \Borel space, $E$ is a countable
  \Borel equivalence relation on $X$, and $\rho \from E \to
  \openinterval{0}{\infty}$ is a \Borel cocycle. Then exactly one of the
  following holds:
  \begin{enumerate}
    \item There is a finite-to-one \Borel compression of $\rho$ over a
      finite \Borel subequivalence relation of $E$.
    \item There is a $\rho$-invariant \Borel probability measure.
  \end{enumerate}
\end{theorem}

\begin{theoremproof}
  Proposition \ref{smooth:compression} ensures that conditions
  (1) and (2) are mutually exclusive. To see that at least one of them
  holds, fix a countable group $\Gamma$ of \Borel automorphisms of
  $X$ whose induced orbit equivalence relation is $E$, and define
  $\rho_\gamma \from X \to \openinterval{0}{\infty}$ by $\rho_\gamma
  (x) = \rho(\gamma \cdot x, x)$ for all $\gamma \in \Gamma$.
  
  By standard change of topology results (see, for example, \cite
  [\S13]{Kechris}), there exist a \Polish topology on
  $\closedopeninterval{0}{\infty}$ and a zero-dimensional \Polish
  topology on $X$, compatible with the underlying \Borel structures
  of $\closedopeninterval{0}{\infty}$ and $X$, with respect to which
  every interval with rational endpoints is clopen, $\Gamma$ acts by
  homeomorphisms, and each $\rho_\gamma$ is continuous. Fix a
  compatible complete metric on $X$, as well as a countable algebra
  $\calU$ of clopen subsets of $X$, closed under multiplication by
  elements of $\Gamma$, and containing a basis for $X$ as well as
  the pullback of every interval with rational endpoints under every
  $\rho_\gamma$.
  
  We say that a function $f \from X \to \closedopeninterval{0}{\infty}$
  is \definedterm{$\calU$-simple} if it is a finite linear combination of
  characteristic functions of sets in $\calU$. Note that for all $\epsilon
  > 0$, $\gamma \in \Gamma$, and $Y \subseteq X$ on which
  $\rho_\gamma$ is bounded, there is such a function with the further
  property that $\absolutevalue{f(y) - \rho_\gamma(y)} \le \epsilon$ for
  all $y \in Y$.
    
  Fix a sequence $\sequence{\epsilon_n}[n \in \N]$ of positive real
  numbers converging to zero, as well as an increasing sequence
  $\sequence{\calU_n}[n \in \N]$ of finite subsets of $\calU$ whose
  union is $\calU$.
  
  By recursively applying Propositions \ref{cocycles:approximation}
  and \ref{cocycles:injection} to functions of the form
  $\equivalenceclass{x}{F} \mapsto \measure{\equivalenceclass{x}
  {F}}{\rho}(A)$ and $\equivalenceclass{x}{F} \mapsto \measure
  {\equivalenceclass{x}{F}}{\rho}(B) - \measure{\equivalenceclass{x}
  {F}}{\rho}(A)$, and throwing out countably-many $E$-invariant
  \Borel sets $B \subseteq X$ for which $\restriction{\rho}{(\restriction
  {E}{B})}$ is smooth, we obtain increasing sequences of finite
  algebras $\calA_n \supseteq \calU_n$ of \Borel subsets of $X$ and
  finite \Borel subequivalence relations $F_n$ of $E$ such that:
  \begin{enumerate}
    \item $\forall n \in \N \forall A \in \calA_n \forall \pair{x}{y} \in E
      \ \measure{\equivalenceclass{x}{F_{n+1}}}{\rho}(A) -
        \measure{\equivalenceclass{y}{F_{n+1}}}{\rho}(A) \le
          \epsilon_n$.
    \item $\forall n \in \N \forall A, B \in \calA_n \ (\forall x \in X
      \ \measure{\equivalenceclass{x}{F_n}}{\rho}(A) \le \measure
        {\equivalenceclass{x}{F_n}}{\rho}(B) \implies$ \\
          \hspace*{5pt} $\exists C \in \calA_{n+1} \forall x \in X \ 0 \le
            \measure{\equivalenceclass{x}{F_{n+1}}}{\rho}(B \setminus
              C)  - \measure{\equivalenceclass{x}{F_{n+1}}}{\rho}(A) \le
                \epsilon_n)$.
  \end{enumerate}
  
  Set $\calA = \union[n \in \N][\calA_n]$ and $F = \union[n \in \N]
  [F_n]$. Condition (1) ensures that we obtain finitely-additive probability
  measures $\mu_x$ on $\calU$ by setting $\mu_x(U) = \lim_{n \goesto
  \infty} \measure{\equivalenceclass{x}{F_n}}{\rho}(U)$ for all $U \in \calU$
  and $x \in X$.
  
  \begin{lemma} \label{invariantprobability:compression:asymptotic}
    Suppose that $\sequence{U_n}[n \in \N]$ is a sequence of
    pairwise disjoint sets in $\calU$ whose union is in $\calU$ and
    $B = \set{x \in X}[{\sum_{n \in \N} \mu_x(U_n) < \mu_x(\union
    [n \in \N][U_n])}]$. Then there is a finite-to-one \Borel compression
    of $\restriction{\rho}{(\restriction{E}{B})}$ over a finite \Borel
    subequivalence relation of $\restriction{E}{B}$.
  \end{lemma}
  
  \begin{lemmaproof}
    As $\mu_x(\union[m \ge n][U_m]) - \sum_{m \ge n} \mu_x(U_m)$
    is independent of $n$, it follows that for all $x \in B$, there exist
    $\delta > 0$ and $n \in \N$ with the property that $\delta + 2
    \sum_{m \ge n} \mu_x(U_m) \le \mu_x(\union[m \ge n][U_m])$. So
    by partitioning $B$ into countably-many $E$-invariant \Borel sets
    and passing to terminal segments of $\sequence{U_n}[n \in \N]$
    on each set, we can assume that $B = \set{x \in X}[{\delta + 2
    \sum_{n \in \N} \mu_x(U_n) \le \mu_x(\union[n \in \N][U_n])}]$ for
    some $\delta > 0$. Fix a sequence $\sequence{\delta_n}[n \in \N]$
    of positive real numbers whose sum is at most $\delta$.
    
    \begin{sublemma}
      There are pairwise disjoint sets $A_n \subseteq \union[m > n]
      [U_m]$ in $\calA$ with the property that for all $n \in \N$, there
      exists $k \in \N$ such that $\forall x \in B \ 0 \le \measure
      {\equivalenceclass{x}{F_k}}{\rho}(A_n) - \measure
      {\equivalenceclass{x}{F_k}}{\rho}(U_n) \le \delta_n$.
    \end{sublemma}
    
    \begin{sublemmaproof}
      Suppose that $n \in \N$ and we have already found $\sequence
      {A_m}[m < n]$. Note that if $x \in B$, then
      \begin{align*}
        \mu_x(U_n) + \sum_{m \ge n} \delta_m
          & \le \mu_x \left( \union[m \in \N][U_m] \right) - \left( \mu_x
            (U_n) + \sum_{m < n} 2\mu_x(U_m) + \delta_m \right) \\
          & \le \mu_x \left( \union[m > n][U_m] \right) - \sum_{m < n}
            \mu_x(U_m) + \delta_m,
      \end{align*}
      so $\forall x \in B \ \measure{\equivalenceclass{x}{F_k}}{\rho}
      (U_n) \le \measure{\equivalenceclass{x}{F_k}}{\rho}(\union[m >
      n][U_m] \setminus \union[m < n][A_m])$ for sufficiently large
      $k \in \N$, by condition (1). It then follows from condition (2) that
      there exists $A_n \subseteq \union[m > n][U_m] \setminus \union
      [m < n][A_m]$ in $\calA$ with the property that $\forall x \in B \ 0
      \le \measure{\equivalenceclass{x}{F_k}}{\rho}(A_n) - \measure
      {\equivalenceclass{x}{F_k}}{\rho}(U_n) \le \delta_n$ for
      sufficiently large $k \in \N$.
    \end{sublemmaproof}
    
    Fix $k_n \in \N$ with the property that $\measure
    {\equivalenceclass{x}{F_{k_n}}}{\rho}(U_n) \le \measure
    {\equivalenceclass{x}{F_{k_n}}}{\rho}(A_n)$ for all $n \in \N$ and
    $x \in B$, as well as \Borel functions $\phi_n \from B \intersection
    U_n \to A_n$ whose graphs are contained in $F_{k_n}$ for all $n
    \in \N$. Then the union of $\union[n \in \N][\phi_n]$ and the
    identity function on $B \setminus \union[n \in \N][U_n]$ is a
    finite-to-one \Borel compression of $\restriction{\rho}{(\restriction
    {E}{B})}$ over the union of $\union[n \in \N][\restriction{F_{k_n}}
    {(A_n \intersection B)}]$ and equality on $B$.
  \end{lemmaproof}
  
  Lemma \ref{invariantprobability:compression:asymptotic} ensures
  that, after throwing out countably-many $E$-invariant \Borel sets
  $B \subseteq X$ for which there is a finite-to-one \Borel compression
  of $\restriction{\rho}{(\restriction{E}{B})}$ over a finite \Borel
  subequivalence relation of $\restriction{E}{B}$, we can assume that
  for all $\delta > 0$ and $U \in \calU$, there is a partition $\sequence
  {U_n}[n \in \N]$ of $U$ into sets in $\calU$ of diameter at most
  $\delta$ such that $\mu_x(U) = \sum_{n \in \N} \mu_x(U_n)$ for all
  $x \in X$.
  
  \begin{lemma}
    Each $\mu_x$ is a measure on $\calU$.
  \end{lemma}
  
  \begin{lemmaproof}
    Suppose, towards a contradiction, that there are pairwise disjoint
    sets $U_n \in \calU$ with $\union[n \in \N][U_n] \in \calU$ but
    $\mu_x(\union[n \in \N][U_n]) > \sum_{n \in \N} \mu_x(U_n)$, for
    some $x \in X$. Fix a sequence $\sequence{\delta_n}[n \in \N]$
    of positive real numbers converging to zero, and recursively
    construct a sequence $\sequence{V_t}[t \in \Bairetree]$ of sets
    in $\calU$, beginning with $V_\emptyset = \union[n \in \N][U_n]$,
    such that $\sequence{V_{t \concatenation \sequence{n}}}[n \in \N]$
    is a partition of $V_t$ into sets of diameter at most
    $\delta_{\length{t}}$ with the property that $\mu_x(V_t) = \sum_{n
    \in \N} \mu_x(V_{t \concatenation \sequence{n}})$, for all $t \in
    \Bairetree$. Set $r = \sum_{n \in \N} \mu_x(U_n)$, and recursively
    construct a sequence $\sequence{i_n}[n \in \N]$ of natural numbers
    with the property that $\sum_{t \in T_n} \mu_x(V_t) > r$, where
    $T_n = \product[m < n][i_m]$, for all $n \in \N$. Set $V_n = \union
    [t \in T_n][V_t]$ for all $n \in \N$. As $\sequence{U_n}[n \in \N]$
    covers the compact set $K = \intersection[n \in \N][V_n]$, so too
    does $\sequence{U_m}[m < n]$, for some $n \in \N$. Set $U =
    \union[m < n][U_m]$, and let $T$ be the tree of all $t \in \union[m
    \in \N][T_m]$ for which $V_t \nsubseteq U$. Note that $T$ is necessarily
    well-founded, since any branch $b$ through $T$ would give rise
    to a singleton $\intersection[n \in \N][V_{\restriction{t}{n}}]$ contained
    in $K \setminus U$. \Koenig's Lemma therefore yields $m \in \N$ with
    $T \subseteq \union[\ell < m][T_\ell]$, in which case $V_m
    \subseteq U$, contradicting the fact that $\mu_x(V_m) > \mu_x
    (U)$.
  \end{lemmaproof}
  
  As a consequence, \Caratheodory's Theorem ensures that
  there is a unique extension of each $\mu_x$ to a \Borel probability
  measure $\completion{\mu}_x$ on $X$.

  \begin{lemma} \label{invariantprobability:compression:rotation}
    Suppose that $\gamma \in \Gamma$, $U \in \calU$,
    $\rho_\gamma$ is bounded on $U$, and $B = \set{x \in X}
    [\completion{\mu}_x(\image{\gamma}{U})
    \neq \int_U \rho_\gamma \ d\completion{\mu}_x]$. Then there is
    a finite-to-one \Borel compression of $\restriction{\rho}
    {(\restriction{E}{B})}$ over a finite \Borel subequivalence relation
    of $\restriction{E}{B}$.
  \end{lemma}
  
  \begin{lemmaproof}
    By the symmetry of our argument, it is enough to establish
    the analogous lemma for the set $B = \set{x \in X}[\completion
    {\mu}_x(\image{\gamma}{U}) < \int_U \rho_\gamma \ d\completion
    {\mu}_x]$. By partitioning $B$ into countably-many $E$-invariant
    \Borel sets, we can assume that $B = \set{x \in X}[\delta +
    \completion{\mu}_x(\image{\gamma}{U}) < \int_U \rho_\gamma
    \ d\completion{\mu}_x]$ for some $\delta > 0$.
    
    \begin{sublemma}
      \label{invariantprobability:compression:rotation:simple}
      For all $\epsilon > 0$, there exists $n \in \N$ with the property
      that $\absolutevalue{\int_U \rho_\gamma \ d\completion{\mu}_x
      - \int_U \rho_\gamma \ d\measure{\equivalenceclass{x}{F_n}}
      {\rho}} \le \epsilon$ for all $x \in X$.
    \end{sublemma}
    
    \begin{sublemmaproof}
      Fix a $\calU$-simple function $f \from X \to \closedopeninterval
      {0}{\infty}$ with the property that $\absolutevalue{f(x) -
      \rho_\gamma(x)} \le \epsilon / 3$ for all $x \in U$. By condition (1),
      there exists $n \in \N$ such that $\absolutevalue{\int_U
      f \ d\completion{\mu}_x - \int_U f \ d\measure{\equivalenceclass{x}
      {F_n}}{\rho}} \le \epsilon / 3$ for all $x \in X$. But then
      \begin{align*}
        \displayedabsolutevalue{\int_U \rho_\gamma \ d\completion
          {\mu}_x - \int_U \rho_\gamma \ d\measure{\equivalenceclass
            {x}{F_n}}{\rho}}
          & \le \displayedabsolutevalue{\int_U \rho_\gamma
            \ d\completion{\mu}_x - \int_U f \ d\completion{\mu}_x} + \\
          & \mspace{21mu} \displayedabsolutevalue{\int_U
              f \ d\completion{\mu}_x - \int_U f \ d\measure
                {\equivalenceclass{x}{F_n}}{\rho}} + \\
          & \mspace{21mu} \displayedabsolutevalue{\int_U f
            \ d\measure{\equivalenceclass{x}{F_n}}{\rho} - \int_U
              \rho_\gamma \ d\measure{\equivalenceclass{x}{F_n}}
                {\rho}} \\
          & \le \epsilon
      \end{align*}
      for all $x \in X$.
    \end{sublemmaproof}
    
    Condition (1) and Sublemma \ref
    {invariantprobability:compression:rotation:simple} ensure that
    there exists $n \in \N$ such that $\measure
    {\equivalenceclass{x}{F_n}}{\rho}(\image{\gamma}{U}) < \int_U
    \rho_\gamma \ d\measure{\equivalenceclass{x}{F_n}}{\rho}$ for
    all $x \in B$. As the former quantity is $\cardinality[x][\rho]{\image
    {\gamma}{U} \intersection \equivalenceclass{x}{F_n}} / \cardinality
    [x][\rho]{\equivalenceclass{x}{F_n}}$ and the latter is $\cardinality
    [x][\rho]{\image{\gamma}{U \intersection \equivalenceclass{x}
    {F_n}}} / \cardinality[x][\rho]{\equivalenceclass{x}{F_n}}$, it follows
    that $\cardinality[x][\rho]{\image{\gamma}{U} \intersection
    \equivalenceclass{x}{F_n}} < \cardinality[x][\rho]{\image{\gamma}
    {U \intersection \equivalenceclass{x}{F_n}}}$ for all $x \in B$, so
    any function from $B \intersection \image{\gamma}{U}$ to $B
    \intersection \image{\gamma}{U}$, sending $\image{\gamma}{U}
    \intersection \equivalenceclass{x}{F_n}$ to $\image{\gamma}{U
    \intersection \equivalenceclass{x}{F_n}}$ for all $x \in B
    \intersection \image{\gamma}{U}$, is a compression of $\restriction
    {\rho}{(\restriction{E}{(B \intersection \image{\gamma}{U})})}$ over
    the equivalence relation $\restriction{\image{(\gamma \times
    \gamma)}{F_n}}{(B \intersection \image{\gamma}{U})}$. The
    \Lusin-\Novikov uniformization theorem yields a \Borel such
    function, and every \Borel such function trivially extends to a
    finite-to-one \Borel compression of $\restriction{\rho}{(\restriction
    {E}{B})}$ over a finite \Borel subequivalence relation of $\restriction
    {E}{B}$.
  \end{lemmaproof}
  
  Lemma \ref{invariantprobability:compression:rotation} ensures that,
  after throwing out countably-many $E$-invariant \Borel sets $B
  \subseteq X$ for which there is a finite-to-one \Borel compression of
  $\restriction{\rho}{(\restriction{E}{B})}$ over a finite \Borel
  subequivalence relation of $\restriction{E}{B}$, we can assume that
  $\completion{\mu}_x(\image{\gamma}{U}) = \int_U \rho_\gamma
  \ d\completion{\mu}_x$ for all $\gamma \in \Gamma$, $U \in
  \calU$ on which $\rho_\gamma$ is bounded, and $x \in X$. As our
  choice of topologies ensures that every open set $U \subseteq X$ is
  a disjoint union of sets in $\calU$ on which $\rho_\gamma$ is
  bounded, we obtain the same conclusion even when $U \subseteq
  X$ is an arbitrary open set. As every \Borel probability measure on a
  \Polish space is regular (see, for example, \cite[Theorem 17.10]
  {Kechris}), we obtain the same conclusion even when $U \subseteq
  X$ is an arbitrary \Borel set. And since every \Borel automorphism
  $T \from X \to X$ whose graph is contained in $E$ is a disjoint
  union of restrictions of automorphisms in $\Gamma$ to \Borel
  subsets, it follows that each $\completion{\mu}_x$ is $\rho$-invariant.
\end{theoremproof}

\begin{acknowledgements}
  I would like to thank the anonymous referee for an unusually
  thorough reading of the original version of the paper and
  many useful suggestions.
\end{acknowledgements}

\bibliographystyle{amsalpha}
\bibliography{bibliography}

\providecommand{\bysame}{\leavevmode\hbox to3em{\hrulefill}\thinspace}
\providecommand{\MR}{\relax\ifhmode\unskip\space\fi MR }
\providecommand{\MRhref}[2]{%
  \href{http://www.ams.org/mathscinet-getitem?mr=#1}{#2}
}
\providecommand{\href}[2]{#2}
\begin{thebibliography}{MVN36}

\bibitem[BK96]{BeckerKechris}
Howard Becker and Alexander~S. Kechris, \emph{The descriptive set theory of
  {P}olish group actions}, London Mathematical Society Lecture Note Series,
  vol. 232, Cambridge University Press, Cambridge, 1996. \MR{1425877}

\bibitem[DJK94]{DoughertyJacksonKechris}
R.~Dougherty, S.~Jackson, and A.~S. Kechris, \emph{The structure of hyperfinite
  {B}orel equivalence relations}, Trans. Amer. Math. Soc. \textbf{341} (1994),
  no.~1, 193--225. \MR{1149121 (94c:03066)}

\bibitem[FM77]{FeldmanMoore}
J.~Feldman and C.~C. Moore, \emph{Ergodic equivalence relations, cohomology,
  and von {N}eumann algebras. {I}}, Trans. Amer. Math. Soc. \textbf{234}
  (1977), no.~2, 289--324. \MR{MR0578656 (58 \#28261a)}

\bibitem[Hop32]{Hopf}
Eberhard Hopf, \emph{Theory of measure and invariant integrals}, Trans. Amer.
  Math. Soc. \textbf{34} (1932), no.~2, 373--393. \MR{1501643}

\bibitem[Kec95]{Kechris}
A.S. Kechris, \emph{Classical descriptive set theory}, Graduate Texts in
  Mathematics, vol. 156, Springer-Verlag, New York, 1995. \MR{1321597
  (96e:03057)}

\bibitem[KM04]{KechrisMiller}
A.~S. Kechris and B.~D. Miller, \emph{Topics in orbit equivalence}, Lecture
  Notes in Mathematics, vol. 1852, Springer-Verlag, Berlin, 2004. \MR{2095154
  (2005f:37010)}

\bibitem[KST99]{KechrisSoleckiTodorcevic}
A.~S. Kechris, S.~Solecki, and S.~Todorcevic, \emph{Borel chromatic numbers},
  Adv. Math. \textbf{141} (1999), no.~1, 1--44. \MR{MR1667145 (2000e:03132)}

\bibitem[Mil08a]{Miller:Probability}
B.~D. Miller, \emph{The existence of measures of a given cocycle. {II}.
  {P}robability measures}, Ergodic Theory Dynam. Systems \textbf{28} (2008),
  no.~5, 1615--1633. \MR{2449547 (2010a:37005)}

\bibitem[Mil08b]{Miller:Sigma}
Benjamin Miller, \emph{The existence of measures of a given cocycle. {I}.
  {A}tomless, ergodic {$\sigma$}-finite measures}, Ergodic Theory Dynam.
  Systems \textbf{28} (2008), no.~5, 1599--1613. \MR{2449546 (2010a:37005)}

\bibitem[MVN36]{MurrayVonNeumann}
F.~J. Murray and J.~Von~Neumann, \emph{On rings of operators}, Ann. of Math.
  (2) \textbf{37} (1936), no.~1, 116--229. \MR{1503275}

\bibitem[Nad90]{Nadkarni}
M.~G. Nadkarni, \emph{On the existence of a finite invariant measure}, Proc.
  Indian Acad. Sci. Math. Sci. \textbf{100} (1990), no.~3, 203--220.
  \MR{1081705}

\bibitem[Wei84]{Weiss}
Benjamin Weiss, \emph{Measurable dynamics}, Conference in modern analysis and
  probability ({N}ew {H}aven, {C}onn., 1982), Contemp. Math., vol.~26, Amer.
  Math. Soc., Providence, RI, 1984, pp.~395--421. \MR{737417 (85j:28027)}

\end{thebibliography}

\end{document}